\documentclass{article}
\usepackage{amsfonts,latexsym}
\usepackage{amsmath, times}
\usepackage{amssymb}
\usepackage{amsthm}

\newcommand{\R}{\mathbb R}

\newcommand{\e}{\varepsilon}
\newtheorem{theorem}{Theorem}[section]
\newtheorem{lemma}{Lemma}
\newtheorem{proposition}{Proposition}[section]

\newtheorem{definition}{Definition}[section]
\newtheorem{remark}{Remark}[section]
\newlength{\defbaselineskip}
\newcommand{\setlinespacing}[2]%
          {\setlength{\baselineskip}{#1 \defbaselineskip}}

\makeatletter

\@addtoreset{equation}{section} \makeatother 
\begin{document}

\begin{center}
 {\Large   {  A Logarithmic Weighted Adams-type inequality in the whole of $\mathbb{R}^{N}$  with an application. }}
\end{center}
\vspace{0.2cm}
\begin{center}
  Rached Jaidane

 \

\noindent\footnotesize  Department of Mathematics, Faculty of Science of Tunis, University of Tunis El Manar, Tunisia.\\
 Address e-mail: rachedjaidane@gmail.com\\
\end{center}

\vspace{0.5cm}
\noindent {\bf Abstract.} In this paper, we will establish a logarithmic weighted Adams inequality in a logarithmic weighted second order Sobolev space in the whole set of $\mathbb{R}^{N}$. Using this result, we delve into the analysis of a weighted fourth-order equation in $\mathbb{R}^{N}$. We assume that the non-linearity of the equation exhibits either critical or subcritical exponential growth, consistent with the Adams-type inequalities previously established. By applying the Mountain Pass Theorem, we demonstrate the existence of a weak solution to this problem. The primary challenge lies in the lack of compactness in the energy caused by the critical exponential growth of the non-linear term $f$.
\\

\noindent {\footnotesize\emph{Keywords:}Adams' inequality, Moser-Trudinger's inequality, $\frac{N}{2}$-biharmonic equation, Nonlinearity of exponential growth, Mountain pass method, Compactness level.\\
\noindent {\bf $2010$ Mathematics Subject classification}: $35$J$20$, $35$J$30$, $35$K$57$, $35$J$60$.}

\section{Introduction}
We start by providing an overview of Trudinger-Moser inequalities within classicalA first-order Sobolev spaces. Additionally, we'll explore Adams' inequalities in second-order Sobolev spaces. Subsequently, we'll extend these concepts to weighted Sobolev spaces.
Moreover, we'll reference relevant works associated with these concepts.\\

In dimension $N\geq 2$ and for bounded domain $\Omega\subset\mathbb{R}^{N}$, the critical exponential growth is given by the well known
Trudinger-Moser inequality \cite{JMo,NST} \begin{align}\label{TM}\displaystyle\sup_{\int_{\Omega} |\nabla u|^{N}\leq1}\int_{\Omega}e^{\alpha|u|^{\frac{N}{N-1}}}dx<+\infty~~\mbox{if and only if}~~\alpha\leq \alpha_{N},\end{align}
where $\alpha_{N}=\omega_{N-1}^{\frac{1}{N-1}}$  with $\omega_{N-1}$ is the area of the unit sphere $S^{N-1}$ in $\mathbb{R}^{N}$.
Later, the Trudinger-Moser inequality was improved to weighted inequalities \cite{ CR1, CR3}.\\
 Equation (\ref{TM}) has been utilized to address elliptic problems that encompass nonlinearities exhibiting exponential growth. For instance, we refer to the following problems in dimensions where $N\geq2$
 \begin{equation}\nonumber
-\Delta_{N} u=-\mbox{div}(|\nabla u|^{N-2}\nabla u)= f(x,u)~~\mbox{in}~~\Omega\subset \mathbb{R}^{N},
\end{equation}
which  have been studied considerably by Adimurthi \cite{Adi1,Adi},  Figueiredo et al. \cite{FMR}, Lam and Lu \cite{LL3, LL2, LL1, LL}, Miyagaki and Souto \cite{MS} and Zhang and Chen \cite{ZC}.\\
Considerable focus has been directed toward weighted inequalities in weighted Sobolev spaces, notably known in mathematical literature as the weighted Trudinger-Moser inequality \cite{CR1, CR3}. The majority of studies have centered on radial functions owing to the radial nature of the weights involved. This quality enhances the maximum growth of integrability.
When the weight is of logarithmic type, Calanchi and Ruf \cite{CR2} extend the Trudinger-Moser inequality and proved the following results in the weighted Sobolev  space,  $W_{0,rad}^{1,N}(B,\rho)=\mbox{closure}\{u \in
C_{0,rad}^{\infty}(B)~~|~~\int_{B}|\nabla u|^{N}\rho(x)dx <\infty\},$ where $B$ denote the unit ball of $\mathbb{R}^{N}$.
\begin{theorem}\cite{CR2} \label{th1.1}\begin{itemize}\item[$(i)$] ~~Let $\beta\in[0,1)$ and let $\rho$ given by $ \rho(x)=\big(\log \frac{1}{|x|}\big)^{\beta}$, then
$$
 \int_{B} e^{|u|^{\gamma}} dx <+\infty, ~~\forall~~u\in W_{0,rad}^{1,N}(B,\rho),~~
  \mbox{if and only if}~~\gamma\leq \gamma_{N,\beta}=\frac{N}{(N-1)(1-\beta)}=\frac{N'}{1-\beta}
$$
and
 $$
 \sup_{\substack{u\in W_{0,rad}^{1,N}(B,\rho) \\ \int_{B}|\nabla u|^{N}w(x)dx\leq 1}}
 \int_{B}~e^{\alpha|u|^{\gamma_{N,\beta} }}dx < +\infty~~~~\Leftrightarrow~~~~ \alpha\leq \alpha_{N,\beta}=N[\omega^{\frac{1}{N-1}}_{N-1}(1-\beta)]^{\frac{1}{1-\beta}}
$$
where $\omega_{N-1}$ is the area of the unit sphere $S^{N-1}$ in $\R^{N}$ and $N'$ is the H$\ddot{o}$lder conjugate of $N$.
\item [$(ii)$] Let $\rho$ given by $\rho(x)=\big(\log \frac{e}{|x|}\big)^{N-1}$, then
  \begin{equation*}\label{eq:71.5}
 \int_{B}exp\{e^{|u|^{\frac{N}{N-1}}}\}dx <+\infty, ~~~~\forall~~u\in W_{0,rad}^{1,N}(B,\rho)
 \end{equation*} and
 \begin{equation*}\label{eq:71.6}
\sup_{\substack{u\in W_{0,rad}^{1,N}(B,\rho) \\  \|u\|_{\rho}\leq 1}}
 \int_{B}exp\{\beta e^{\omega_{N-1}^{\frac{1}{N-1}}|u|^{\frac{N}{N-1}}}\}dx < +\infty~~~~\Leftrightarrow~~~~ \beta\leq N,
 \end{equation*}
where $\omega_{N-1}$ is the area of the unit sphere $S^{N-1}$ in $\R^{N}$ and $N'$ is the H$\ddot{o}$lder conjugate of $N$.\end{itemize}
\end{theorem}
The theorem (\ref{th1.1}) has enabled the exploration of second-order weighted elliptic problems in dimensions where $N\geq2$. As a result, Calanchi et al.\cite{CRS} established the existence of a non-trivial radial solution for an elliptic problem defined on the unit ball in $\mathbb{R}^2$, where the nonlinearities exhibit double exponential growth at infinity. Following this, Deng et al. investigated the subsequent problem \begin{equation}\label{pro}
\displaystyle \left\{
\begin{array}{rclll}
-\textmd{div} (\sigma(x)|\nabla u(x)|^{N-2}\nabla u(x) ) &=& \ f(x,u)& \mbox{in} & B \\
u&=&0 &\mbox{on }&  \partial B,
\end{array}
\right.
 \end{equation}
  where $B$ is the unit ball in $\R^N,\; N\geq2$ and the nonlinearity $f(x, u)$ is continuous in $B\times\R$ and has critical growth in the sense of Theorem (\ref{th1.1}).The authors have proved that there is a non-trivial solution to this problem, using the mountain pass Theorem. Similar results are proven by Chetouane and Jaidane \cite{CJ}, Dridi \cite{DRI1} and Zhang \cite{Z}. Furthermore, problem \eqref{pro}, involving a potential, has been studied by Baraket and Jaidane \cite{BJ}.
Moreover, we mention that Abid et al. \cite{ABJ} have proved the existence of a positive ground state solution for a weighted second-order elliptic problem of Kirchhoff type, with nonlinearities having a double exponential growth at infinity, using minimax techniques combined with Trudinger-Moser inequality.\\
 \\
In recent years, Aouaoui and Jlel \cite{Aou} have extended the work of Calanchi and Ruf to the whole $\mathbb{R}^2$ space, by considering the following weight
\begin{equation} \label{1.1}
  \rho_{\beta}(x)= \begin{cases}\left(\log \left(\frac{e}{|x|}\right)\right)^{\beta} & \text { if }|x|<1, \\ \chi(|x|) & \text { if }|x| \geq 1,\end{cases}
\end{equation}

where, $0<\beta \leq 1$ and $\chi:[1,+\infty[\rightarrow] 0,+\infty[$ is a continuous function such that $\chi(1)=1$ and $\inf _{t \in[1,+\infty[} \chi(t)>0$. The authors consider the space $E_{\beta}$ as the space of all radial functions of the completion of $C_{0}^{\infty}\left(\mathbb{R}^{2}\right)$ with respect to the norm

$$
\|u\|_{E_{\beta}}^{2}=\int_{\mathbb{R}^{2}}|\nabla u|^{2} \rho_{\beta}(x) d x+\int_{\mathbb{R}^{2}} u^{2} d x=|\nabla u|_{L^{2}\left(\mathbb{R}^{2}, \rho_{\beta}\right)}^{2}+|u|_{L^{2}\left(\mathbb{R}^{2}\right)}^{2} .$$
The authors proved the following result:
\begin{theorem} Let $\beta \in(0,1)$ and $\rho_{\beta}$ be defined by \eqref{1.1}. For all $u \in E_{\beta}$, we have

$$
\int_{\mathbb{R}^{2}}\left(e^{|u|^{\frac{2}{1-\beta}}}-1\right) d x<+\infty.
$$
Moreover, if $\alpha<\tau_{\beta}$, then
\begin{equation} \label{1.5}
    \sup _{u \in E_{\beta},\|u\|_{E_{\beta}} \leq 1} \int_{\mathbb{R}^{2}}\left(e^{\alpha|u|^{\frac{2}{1-\beta}}}-1\right) d x<+\infty
\end{equation}

where $\tau_{\beta}=2\big[2\pi (1-\beta)\big]\frac{1}{1-\beta}$.\\
If $\alpha>\tau_{\beta}$, then

$$
\sup _{u \in E_{\beta},\|u\|_{E_{\beta}} \leq 1} \int_{\mathbb{R}^{2}}\left(e^{\alpha|u|^{\frac{2}{1-\beta}}}-1\right) d x=+\infty .
$$
\end{theorem}
The concept of critical exponential growth was further expanded to higher order Sobolev spaces by Adams \cite{Ada}. Specifically, Adams demonstrated the following outcome: for $m \in \mathbb{N}$ and $\Omega$ as an open bounded set in $\mathbb{R}^N$ where $m<N$, there exists a positive constant $C_{m,N}$ such that\begin{equation}\label{eq:1.31}
\displaystyle\sup_{{u\in W_{0}^{m,\frac{N}{m}}(\Omega), |\nabla^m u|_{\frac{N}{m}} \leq 1}}
\int_{\Omega}~e^{\displaystyle \beta_{0} |u|^{\frac{N}{N-m}}}~dx \leq C_{m,N} \vert \Omega\vert,
\end{equation}
where   $ W_{0}^{m,\frac{N}{m}}(\Omega)$  denotes the $m^{th}$-order Sobolev space, $\nabla^m u$ denotes the $m^{th}$-order gradient of $u$, namely
$$\nabla^m u := \left\{ \begin{array}{ll} \displaystyle\Delta^{\frac{m}{2}} u, \qquad \mbox{if}\; m\; \mbox{is even} \\\\\displaystyle
\nabla\Delta^{\frac{m-1}{2}} u , \qquad \mbox{if}\; m \;\mbox{odd}
\end{array}\right. $$
and
$$\beta_{0} = \beta_{0}(m,N) :=  \frac{N}{\omega_{N-1}} \left\{ \begin{array}{ll}\displaystyle\Big[ \frac{\pi^{\frac{N}{2}} 2^m \Gamma (\frac{m}{2}) }{\Gamma (\frac{N-m}{2}) }\Big]^{\frac{N}{N-m}}, \qquad \mbox{if}\; m\; \mbox{is even} \\\\ \displaystyle
\Big[ \frac{\pi^{\frac{N}{2}} 2^m \Gamma (\frac{m+1}{2}) }{\Gamma (\frac{N-m+1}{2}) }\Big]^{\frac{N}{N-m}}, \qquad \mbox{if}\; m \;\mbox{odd}.
\end{array}\right. $$

  In the particular case where $N = 4$ and
$m = 2,$ the inequality \eqref{eq:1.31} takes the form
\begin{equation}\label{eq:1.311}
\displaystyle\sup_{{u\in W_{0}^{2,2}(\Omega), |\Delta u|_{2} \leq 1}}
\int_{\Omega}~e^{\displaystyle 32 \pi^2 u^{2}}~dx \leq C \vert \Omega\vert.
\end{equation}

 Also, for bounded domains $\Omega\subset \mathbb{R}^{4}$, in \cite{Ada,BRFS} the authors proved the following inequality
 \begin{equation*}
 \sup_{\substack{u\in S }}
 \int_{\Omega}~(e^{\alpha u^{2}}-1)dx < +\infty~~~~\Leftrightarrow~~~~ \alpha\leq 32 \pi^{2}
 \end{equation*}
 where $$S=\{u\in W^{2,2}_{0}(\Omega)~~|~~\displaystyle\big(\int_{\Omega}|\Delta u|^{2}dx \big)^{\frac{1}{2}}\leq 1\}.$$
 When $\Omega$
 is replaced by the whole space $\R^{4}$, Ruff and Sani  [38]
established the corresponding Adams type inequality as follows:\begin{equation}\label{RS}
 \sup_{\substack{\|u \|_{W^{2,2}}}\leq1}
 \int_{\Omega}~(e^{\alpha u^{2}}-1)dx < +\infty~~~~\Leftrightarrow~~~~ \alpha\leq 32 \pi^{2}
 \end{equation}
 where $ \displaystyle\|u\|^{2}_{W^{2,2}(\R^{4})}=\int_{\mathbb{R}^{4}}|\triangle u|^{2}  d x+\int_{\mathbb{R}^{4}}|\nabla u|^{2}dx+\int_{\mathbb{R}^{4}} u^{2} d x $.\\

 Recently,  Adams’ type inequalities on the logarithmic weighted Sobolev space $W^{2,2}_{0,rad}(B_{1})$ of radial function in the unit ball $B$ of $\R^{4}$ has been established. More precisely, in \cite{WZ} the authors proved the following result\begin{theorem}\cite{WZ} \label{th1.11} ~~Let $\beta\in(0,1)$ and let $w=(\log(\frac{e}{|x|}))^{\beta}$, then

 \begin{equation*}
 \sup_{\substack{u\in W_{0,rad}^{2,2}(B_{1},w) \\  \int_{B_{1}}w(x)|\Delta u|^{2}dx \leq 1}}
 \int_{B_{1}}~e^{\displaystyle\alpha|u|^{\frac{2}{1-\beta} }}dx < +\infty~~~~\Leftrightarrow~~~~ \alpha\leq \alpha_{\beta}=4[8\pi^{2}(1-\beta)]^{\frac{1}{1-\beta}}
 \end{equation*}

 \end{theorem}
 This last result allowed the authors in \cite{DRI2} to study the following weighted problem \begin{equation*}\label{eq:1.1}
  \displaystyle \left\{
      \begin{array}{rclll}
  \Delta(w(x) \Delta u)-\Delta u+V(x)u &=&  \displaystyle f(x,u)& \mbox{in} & B_{1} \\
        u=\frac{\partial u}{\partial n}&=&0 &\mbox{on }&  \partial B_{1},
      \end{array}
    \right.
\end{equation*}
where the weight $w(x)$ is given by \begin{equation*}\label{eq: 1.2}
w(x)=(\log\frac{e}{\vert x\vert})^{\beta},  \beta\in (0,1),
\end{equation*}$B_{1}$ is the unitary disk in $\mathbb{R}^{4}$, $f(x,t)$ is continuous in
$B_{1}\times \mathbb{R}$ and behaves like $\exp\{\alpha t^{\frac{2}{1-\beta}}\}$ as $ |t|\rightarrow+ \infty$, for some $\alpha >0$ uniformly with respect to $x\in B_{1}$.   The potential $V :\overline{B_{1}}\rightarrow \mathbb{R}$ is a positive continuous function and bounded away from zero in $B_{1}$.
  The authors establish the existence of radial solution by variational techniques and using Adams' inequality \cite{WZ}. \\
  It should be noted that several works concerning weighted elliptic equations of Kirchhof type with critical nonlinearities in the sense of Theorem \ref{th1.1}  or Theorem \ref{th1.11} have been studied (see \cite{ABJ,RJT,J}).\\

  Denote by $E$  as the space of all radial functions of the completion of $C^{\infty}_{0}(\mathbb{R}^{N})$ with respect to the norm
$$
\|u\|^{\frac{N}{2}}=\int_{\mathbb{R}^{N}}|\triangle u|^{\frac{N}{2}} w_{\beta}(x) d x+\int_{\mathbb{R}^{N}}|\nabla u|^{\frac{N}{2}}dx+\int_{\mathbb{R}^{N}} |u|^{\frac{N}{2}} d x .
$$ where the weight $w_{\beta}(x)$ is given by \begin{equation} \label{w}
  w_{\beta}(x)= \begin{cases}\left(\log \left(\frac{e}{|x|}\right)\right)^{\beta(\frac{N}{2}-1)} & \text { if }|x|<1, \\  \chi(|x|) & \text { if }|x| \geq 1,\end{cases}
\end{equation} where, $\frac{N^{2}-4N+2}{N(N-2)}<\beta < 1$ and $\chi:[1,+\infty[\rightarrow[ 1,+\infty[$ is a continuous function such that $\chi(1)=1$ and $\inf _{t \in[1,+\infty[} \chi(t)\geq1$. Also, we suppose that there exists a positive constant  $M>0$ such that \begin{equation} \label{x1}
    \frac{1}{r^{\frac{N^{2}}{2}}}\left(\int_{1}^{r} t^{N-1} \chi(t) d t\right)\left(\int_{1}^{r} \frac{t^{N-1}}{\chi(t)} d t\right)^{\frac{N}{2}-1} \leq M, \forall r \geq 1,
\end{equation}

\begin{equation} \label{x2}
    \frac{1}{r^{\frac{N^{2}}{2}}}\left(\int_{1}^{r} t^{N-1} \chi(t) d t\right) \leq M, \forall r \geq 1,
\end{equation}

and

\begin{equation} \label{x3}
    \frac{\max_{r \leq t \leq 4 r} \chi(t)}{\min_{r \leq t \leq 4 r} \chi(t)} \leq M, \forall r \geq 1,
\end{equation}
We give some examples of functions $\chi:[1,+\infty[\rightarrow[ 1,+\infty[$ satisfying the conditions (\ref{x1}), (\ref{x2}) and (\ref{x3}):\\
$\bullet$ Any function continuous $\chi$ such that $\chi(1)=1$ and $$1\leq \displaystyle\inf_{t\geq1}\chi(t)\leq \displaystyle \sup_{t\geq1}\chi(t)<+\infty.$$
$\bullet$ $\chi(t)=t^{\delta}, 0<\delta<\frac{N^{2}}{2}-N.$\\
$\bullet$ $\chi(t)=1+\log^{\sigma} t~~, \sigma>1$.\\
Since  the weight $w_{\beta}$ belongs to the Muckenhoupt's class $A_{\frac{N}{2}}$, then $C_{0}^{\infty}\left(\mathbb{R}^{N}\right)$ is dense in the space $E$ (see Lemma \ref{a2}). It follows that the space $E$ can be seen as $$E=
\left\{u \in L_{rad}^{\frac{N}{2}}\left(\mathbb{R}^{N}\right), \int_{\mathbb{R}^{N}}\big(|\triangle u|^{\frac{N}{2}}w_{\beta}(x)+~| \nabla u|^{\frac{N}{2}}\big) d x<+\infty\right\}
,$$ endowed with the norm$$
\|u\|^{\frac{N}{2}}=\int_{\mathbb{R}^{N}}|\triangle u|^{\frac{N}{2}} w_{\beta}(x) d x+\int_{\mathbb{R}^{N}}|\nabla u|^{\frac{N}{2}}dx+\int_{\mathbb{R}^{N}} |u|^{\frac{N}{2}} d x .
$$

In this paper, we prove a weighted Adams' inequality analogous to (\ref{RS}) in the whole of $ \mathbb{R}^{4}$ that is:
  \begin{theorem}\label{th1.4}~~Let $\beta\in(\frac{N^{2}-4N+2}{N(N-2)},1)$ and let $w$ given by (\ref{w}), then

\begin{itemize}
\item[(i)]\begin{equation}\label{subc}\displaystyle\int_{\mathbb{R}^{4}}\big(e^{|u|^{\gamma}}-1\big)dx<+\infty,~~\forall u\in E.\end{equation}
\item[(ii)] \begin{equation}\label{eq:1.3}
 \sup_{\substack{u\in W_{0,rad}^{2,\frac{N}{2}}(B,w) \\  \int_{B}w(x)|\Delta u|^{\frac{N}{2}}dx \leq 1}}
 \int_{B}~(e^{\displaystyle\alpha|u|^{\frac{N}{(N-2)(1-\beta)} }}-1)dx < +\infty~~~~\Leftrightarrow~~~~ \alpha\leq \alpha_{\beta}
\end{equation}
\end{itemize}
$\mbox{ with}~~ \alpha_{\beta}=N[(N-2)N V_{N}]^{\frac{2}{(N-2)(1-\beta)}}(1-\beta)^{\frac{1}{(1-\beta)}}\cdot$
   \end{theorem}

As an application of this last result, we study the following weighted problem
\begin{equation}\label{eq:1.1}
  \displaystyle
  L:=\Delta(w(x)|\Delta u|^{\frac{N}{2}-2} \Delta u)-\text{div}(|\nabla u|^{\frac{N}{2}-2}\nabla u)+V(x)|u|^{\frac{N}{2}-2}u= \displaystyle f(u)~~\mbox{in} ~~ \mathbb{R}^{N},
      \end{equation}
where the weight is given by (\ref{w}), $N\geq4$.The non linearity  $f(t)$ is continuous in
$\mathbb{R}$ and behaves like $\exp\{\alpha t^{\frac{N}{(N-2)(1-\beta)}}\}$ as $ |t|\rightarrow+ \infty$, for some $\alpha >0$ . \\
 Let $\gamma:=\displaystyle\frac{N}{(N-2)(1-\beta)}$. In view  of inequality (\ref{eq:1.3}), we say that $f$ has critical growth at $+\infty$ if there exists some $\alpha_{0}>0$,
\begin{equation}\label{eq:1.4}
\lim_{|s|\rightarrow +\infty}\frac{|f(s)|}{e^{\alpha s^{\gamma}}}=0,~~~\forall~\alpha~~\mbox{such that}~~ \alpha>\alpha_{0} ~~~~
\mbox{and}~~~~\lim_{|s|\rightarrow +\infty}\frac{|f(s)|}{e^{\alpha s^{\gamma}}}=+\infty,~~\forall~ \alpha<\alpha_{0} .\\
\end{equation}
 In view  of inequality (\ref{subc}), we say that $f$ has subcritical growth at $+\infty$ if
 $$\lim_{|s|\rightarrow +\infty}\frac{|f(s)|}{e^{\alpha s^{\gamma}}}=0,~~~\forall~\alpha>0.$$

Let us now state our results. In this paper, we always assume that  the nonlinearities $f(t) $ has critical  growth with $\alpha_{0}> 0$ or  $f(t) $ has subcritical growth and satisfies these conditions:
\begin{description}
\item[$(H_{1})$] The non-linearity $f:  \mathbb{R}\rightarrow\mathbb{R}$ is continuous.
  \item[$(H_{2})$]$\mbox{There exists } \theta > N,  \mbox{ such that } 0 <\theta F( t)
=\displaystyle\theta \int_{0}^{t}f(s)ds \leq tf(t),~~ \forall t\in \mathbb{R}\setminus\{0\}.$
    \item[$(H_{3})$] $\displaystyle\lim_{t\rightarrow 0}\frac{f(t)}{t^{\frac{N}{2}-1}}=0.$

\item[$(H_{4})$] There exist $t_{0},~M_{0}>0$ such that
\begin{align}\nonumber 0<F(t)\leq M_{0} |f(t)|~~\mbox{for all}~~|t|\geq t_{0}.\end{align}
\item[$(H_{5})$] The asymptotic condition \begin{align}\nonumber \displaystyle\lim_{t\rightarrow \infty}\frac{f(t)t}{e^{\alpha_{0}t^{\gamma}}}\geq
\gamma_{0} ~~\mbox{with}~~\gamma_{0}>\displaystyle \frac{(\frac{\alpha_{\beta}}{\alpha_{0}})^{\frac{N}{2\gamma}}}{V_{N}e^{N(1-\log (2e))}}.\end{align}

\end{description}

We say that $u$ is a solution to the problem (\ref{eq:1.1}), if $u$ is a weak solution
in the following sense.
\begin{definition}\label{def1} A function $u$ is called a solution to $ (\ref{eq:1.1}) $ if $u \in E$ and

\begin{equation}\label{eq:1.6}
\int_{\mathbb{R}^{N}}\big(w(x)~|\Delta u|^{\frac{N}{2}-2}~\Delta u~~\Delta\varphi+|\nabla u|^{\frac{N}{2}-2} \nabla u.\nabla \varphi +|u|^{\frac{N}{2}-2}u\varphi\big)~dx=\int_{\mathbb{R}^{N}}f(u)~ \varphi~dx,~~~~~~\mbox{for all }~~\varphi \in E.
\end{equation}
\end{definition}
It is easy to see that seeking weak solutions of the problem (\ref{eq:1.1}) is
 equivalent to find nonzero critical points of the following functional on $ E$:
 \begin{equation}\label{energy}
\mathcal{J}(u)=\frac{2}{N}\bigg(\int_{\mathbb{R}^{N}} w_{\beta}(x) |\Delta u|^{\frac{N}{2}} +|\nabla u|^{\frac{N}{2}}+|u|^{\frac{N}{2}}dx\bigg) -\int_{\mathbb{R}^{N}}F(x,u)dx,
 \end{equation}
 where $F(u)=\displaystyle\int_{0}^{u}f(t)dt$.\\

We prove the following results:\\

Also, in the critical (respectively subcritical) case  we prove the following theorems:
\begin{theorem}\label{th1.5}
Assume that the function $f$  has a critical growth at
$+\infty$ and satisfies
the conditions $(H_{1})$, $(H_{2})$, $(H_{3})$, $(H_{4})$ and  $(H_{5})$ .
Then the problem  (\ref{eq:1.1}) has a nontrivial solution.
\end{theorem}
\begin{theorem}\label{th1.6}
Assume that the function $f$  has subcritical growth at
$+\infty$ and satisfies
the conditions $(H_{1})$, $(H_{2})$, $(H_{3})$, and $(H_{4})$  .
Then the problem  (\ref{eq:1.1}) has a nontrivial solution.
\end{theorem}
In general the study of  fourth order partial differential equations is considered an interesting topic. The interest in studying such equations was stimulated by their applications in micro-electro-mechanical systems, phase field models of multi-phase systems, thin film theory, surface diffusion on solids, interface dynamics, flow in Hele-Shaw cells, see \cite{D, FW, M}. However  many applications are generated
by elliptic problems, such as  the study of traveling waves in suspension bridges, radar imaging  (see, for example \cite{AEG, LL}).

 This paper is structured as follows:

Section 2 presents essential background information on functional spaces.
Section 3 establishes preliminary results crucial for our proofs.
Section 4 focuses on proving Theorem \ref{th1.4}.
Section 5 demonstrates a concentration compactness result akin to Lions' theorem.
Section 6 verifies that the energy $\mathcal{J}$ adheres to two specific geometric properties and a compactness condition, albeit under a specified level.
Section 7 offers the proof of Theorem \ref{th1.4}.
Finally, in Section 8, we conclude by proving Theorem \ref{th1.5} and Theorem \ref{th1.6}.
Through this paper, the constants $C$ or $c$ may change from line to another and we sometimes index the constants in order to show how they change.
\section{Weighted Lebesgue and Sobolev Spaces setting}
Let $\Omega \subset \R^{N}$, $N\geq2$, bounded or unbounded, possibly even equal to the whole  $\R^{N}$  and let $w\in L^{1}(\Omega)$ be a nonnegative function. In order to work with a weighted operator, it becomes necessary to introduce specific functional spaces denoted as $L^{p}(\Omega,w)$, $W^{m,p}(\Omega,w)$, and $W_{0}^{m,p}(\Omega,w)$. Later on, these spaces and some of their properties will be utilized. Let $S(\Omega)$ be the set of all measurable real-valued functions defined on $\Omega$ and two measurable functions are considered as the same element if they are equal almost everywhere.
Following  Drabek et al. \cite{DKN} and Kufner in \cite{Kuf}, the weighted Lebesgue space $L^{p}(\Omega,w)$ is defined as follows:
$$L^{p}(\Omega,w)=\{u:\Omega\rightarrow \R ~\mbox{measurable};~~\int_{\Omega} w_{\beta}(x)|u|^{p}~dx<\infty\}$$
for any real number $1\leq p<\infty$.\\
This is a normed vector space equipped with the norm
$$\|u\|_{p,w}=\Big(\int_{\Omega}w(x)|u|^{p}~dx\Big)^{\frac{1}{p}}.$$
For $m\geq 2$, let $w$ be a given family of weight functions $w_{\tau}, ~~|\tau|\leq m,$ $w=\{w_{\tau}(x)~~x\in\Omega,~~|\tau|\leq m\}.$\\

In \cite{DKN}, the  corresponding weighted Sobolev space was  defined as
$$ W^{m,p}(\Omega,w)=\{ u \in L^{p}(\Omega), D^{\tau} u \in L^{p} (\Omega,w)~~  \forall ~~1\leq|\tau|\leq m-1 , D^{\tau} u   \in L^{p}(\Omega,w) ~~  \forall ~~|\tau|=m  \}$$
endowed with the following norm:

\begin{equation*}\label{eq:2.2}
\|u\|_{W^{m,p}(\Omega,w)}=\bigg(\sum_{ |\tau|\leq m-1}\int_{\Omega}|D^{\tau}u|^{p}dx+\displaystyle \sum_{ |\tau|= m}\int_{\Omega}w(x) |D^{\tau}u|^{p}dx\bigg)^{\frac{1}{p}},
\end{equation*}
where $w_{\tau}=1~~\mbox{for all}~~|\tau|< k,$ $w_{\tau}=w~~\mbox{for all}~~|\tau|=k$.\\
If we suppose also that $w(x)\in L^{1}_{loc}(\Omega)$, then $C^{\infty}_{0}(\Omega)$ is a subset of $W^{m,p}(\Omega,w)$ and we can introduce the space $$W^{m,p}_{0}(\Omega,w)$$
as the closure of $C^{\infty}_{0}(\Omega)$ in $W^{m,p}(\Omega,w).$  Moroever, the injection $$W^{m,p}(\Omega,w)\hookrightarrow W^{m-1,p}(\Omega)~~\mbox{is compact}.$$
Also, $(L^{p}(\Omega,w),\|\cdot\|_{p,w})$ and $(W^{m,p}(\Omega,w),\|\cdot\|_{W^{m,p}(\Omega,w)})$ are separable, reflexive Banach spaces provided that $w(x)^{\frac{-1}{p-1}} \in L^{1}_{loc}(\Omega)$.
Then the space $$E=\{u\in  L^{\frac{N}{2}}_{rad}(\mathbb{R}^{N})~~|~~\int_{\mathbb{R}^{N}}\big(w_{\beta}(x)|\Delta u|^{\frac{N}{2}}+| \nabla u|^{\frac{N}{2}}\big)dx<+\infty \}$$ is a Banach and reflexive space. \\
We have the following result \begin{lemma}\label{a2}  $C_{0}^{\infty}\left(\mathbb{R}^{N}\right)$ is dense in the space

$$
\left\{u \in L^{\frac{N}{2}}\left(\mathbb{R}^{N}\right), \int_{\mathbb{R}^{N}}\big(|\triangle u|^{\frac{N}{2}}w_{\beta}(x)+~| \nabla u|^{\frac{N}{2}}\big) d x<+\infty\right\}
$$\end{lemma}
\textit{Proof} it suffice to see that $\omega_{\beta}$ belongs to the Muckenhoupt's class $A_{\frac{N}{2}}$ (we also say that $\omega_{\beta}$ is an $A_{\frac{N}{2}}$-weight), that is

$$
\sup \left(\frac{1}{|B|} \int_{B} w_{\beta}(x) d x\right)\left(\frac{1}{|B|} \int_{B}\left(w_{\beta}(x)\right)^{-1} d x\right)^{\frac{N}{2}-1}<+\infty
,$$ where the supremum is taken over all balls $B \subset \mathbb{R}^{N}$.\\

Let $r>0$ and $x_{0} \in \mathbb{R}^{4}$. Denote by $B\left(x_{0}, r\right)$ (resp. $\left.B(0, r)\right)$ the open ball of $\mathbb{R}^{4}$ of center $x_{0}$ and radius $r$ (resp. of center 0 and radius $r$ ).

First case: Suppose that $B\left(x_{0}, r\right) \cap B(0, r) \neq \emptyset$. Thus, $B\left(x_{0}, r\right) \subset B(0,3 r)$ which implies that

$$
\frac{1}{\left|B\left(x_{0}, r\right)\right|^{\frac{N}{2}}}\left(\int_{B\left(x_{0}, r\right)} w_{\beta}(x) d x\right)\left(\int_{B\left(x_{0}, r\right)} \frac{d x}{w_{\beta}(x)}\right)^{\frac{N}{2}-1} $$
\begin{equation} \label{7.1}
    \leq  \frac{c}{r^{\frac{N^{2}}{2}}}\left(\int_{0}^{3 r} w_{\beta}(t) t^{N-1} d t\right)\left(\int_{0}^{3 r} \frac{t^{N-1}}{w_{\beta}(t)} d t\right)^{\frac{N}{2}-1} .
\end{equation}

If $3 r<1$, then

$$
\begin{aligned}
& \frac{c}{r^{\frac{N^{2}}{2}}}\left(\int_{0}^{3 r} w_{\beta}(t) t^{N-1} d t\right)\left(\int_{0}^{3 r} \frac{t^{N-1}}{w_{\beta}(t)} d t\right)^{\frac{N}{2}-1} \\
= & \frac{c}{r^{\frac{N^{2}}{2}}}\left(\int_{0}^{3 r} t^{N-1}(1-\log t)^{\beta(\frac{N}{2}-1)} d t\right)\left(\int_{0}^{3 r} \frac{t^{N-1}}{(1-\log t)^{\beta(\frac{N}{2}-1)}} d t\right)^{\frac{N}{2}-1} .
\end{aligned}
$$

But, a simple computation gives

\begin{equation} \label{7.2}
  \limsup _{r \rightarrow 0^{+}} \frac{c}{r^{\frac{N^{2}}{2}}}\left(\int_{0}^{3 r} t^{N-1}(1-\log t)^{\beta(\frac{N}{2}-1)} d t\right)\left(\int_{0}^{3 r} \frac{t^{N-1}}{(1-\log t)^{\beta(\frac{N}{2}-1)}} d t\right)^{\frac{N}{2}-1}<+\infty .
\end{equation}

If $3 r \geq 1$, then

$$
 \frac{c}{r^{\frac{N^{2}}{2}}}\left(\int_{0}^{3 r} \omega_{\beta}(t) t^{N-1} d t\right)\left(\int_{0}^{3 r} \frac{t^{N-1}}{w_{\beta}(t)} d t\right)^{\frac{N}{2}-1}
$$

\begin{equation} \label{7.3}
    = \frac{c}{r^{\frac{N^{2}}{2}}}\left(\int_{0}^{1} t^{N-1}(1-\log t)^{\beta(\frac{N}{2}-1)} d t+\int_{1}^{3 r} t \chi(t) d t\right)\left(\int_{0}^{1} \frac{t^{N-1}}{(1-\log t)^{\beta(\frac{N}{2}-1)}} d t+\int_{1}^{3 r} \frac{t^{N-1}}{\chi(t)} d t\right)^{\frac{N}{2}-1} .
\end{equation}

Since $\inf _{t \geq 1} \chi(t)\geq1$, then

\begin{equation}\label{7.4}
    \limsup _{r \rightarrow+\infty} \frac{1}{r^{\frac{N^{2}}{2}}} \int_{1}^{3 r} \frac{t^{N-1}}{\chi(t)} d t=0<+\infty
\end{equation}

On the other hand, by \eqref{x2}, we infer

\begin{equation}\label{7.5}
    \limsup _{r \rightarrow+\infty} \frac{1}{r^{\frac{N^{2}}{2}}} \int_{1}^{3 r} t^{N-1} \chi(t) d t<+\infty
\end{equation}

Hence, in view of \eqref{7.4}, \eqref{7.5} and \eqref{7.3}, it remains to show that

$$
\limsup _{r \rightarrow+\infty} \frac{1}{r^{\frac{N^{2}}{2}}}\left(\int_{1}^{3 r} t^{N-1} \chi(t) d t\right)\left(\int_{1}^{3 r} \frac{t^{N-1}}{\chi(t)} d t\right)^{\frac{N}{2}-1}<+\infty .
$$

But this fact can immediately be deduced from \eqref{x2}. Combining \eqref{7.2} and \eqref{7.3}, we deduce from \eqref{7.1} that there exists a constant $D_{0}>0$ independent of $x_{0}$ and $r$ such that

\begin{equation}\label{7.6}
    \frac{1}{\left|B\left(x_{0}, r\right)\right|^{\frac{N}{2}}}\left(\int_{B\left(x_{0}, r\right)} \omega_{\beta}(x) d x\right)\left(\int_{B\left(x_{0}, r\right)} \frac{d x}{\omega_{\beta}(x)}\right)^{\frac{N}{2}-1} \leq D_{0}
\end{equation}

Second case: Suppose that $B\left(x_{0}, r\right) \cap B(0, r)=\emptyset$. In this case, we have

$$
\frac{\left|x_{0}\right|}{2} \leq|x| \leq 2\left|x_{0}\right|, \forall x \in B\left(x_{0}, r\right) \text {. }
$$

Hence,

$$
 \frac{1}{\left|B\left(x_{0}, r\right)\right|^{\frac{N}{2}}}\left(\int_{B\left(x_{0}, r\right)} w_{\beta}(x) d x\right)\left(\int_{B\left(x_{0}, r\right)} \frac{d x}{w_{\beta}(x)}\right)^{\frac{N}{2}-1}$$
\begin{equation}\label{7.7}
   \displaystyle \leq  \left(    \frac{\sup_{\frac{|x_0|}{2} \leq |x| \leq  2 |x_0|} w_{\beta}(t)}{\inf_{\frac{|x_0|}{2} \leq |x| \leq  2 |x_0|} w_{\beta}(t)} \right)     \leq    \sup_{\tau > 0} \left(    \frac{\sup_{\tau \leq t \leq 4 \tau} w_{\beta}(t)}{\inf_{\tau \leq t \leq 4 \tau} w_{\beta}(t)} \right)
\end{equation}

If $4 \tau<1$, then

$$
\frac{\sup _{\tau \leq t \leq 4 \tau} w_{\beta}(t)}{\inf _{\tau \leq t \leq 4 \tau}w_{\beta}(t)}=\frac{(1-\log \tau)^{\beta(\frac{N}{2}-1)}}{(1-\log (4 \tau))^{\beta(\frac{N}{2}-1)}} .
$$

Taking into account that

$$
\sup _{0<\tau<\frac{1}{4}}\left(\frac{1-\log \tau}{1-\log (4 \tau)}\right)^{\beta(\frac{N}{2}-1)}<+\infty
$$

it follows that

\begin{equation} \label{7.8}
\sup_{0<\tau<\frac{1}{4}}\left(\frac{\sup_{\tau \leq t \leq 4 \tau} w_{\beta}(t)}{\inf_{\tau \leq t \leq 4 \tau}w_{\beta}(t)}\right)<+\infty .
\end{equation}

If $\frac{1}{4} \leq \tau<1$, then

$$
\frac{\sup _{\tau \leq t \leq 4 \tau}w_{\beta}(t)}{\inf _{\tau \leq t \leq 4 \tau}w_{\beta}(t)} \leq \frac{\sup _{\frac{1}{4} \leq t \leq 4}w_{\beta}(t)}{\inf _{\frac{1}{4} \leq t \leq 4}w_{\beta}(t)}<+\infty,
$$

and consequently,

\begin{equation} \label{7.9}
    \sup _{\frac{1}{4} \leq \tau<1}\left(\frac{\sup _{\tau \leq t \leq 4 \tau}w_{\beta}(t)}{\inf _{\tau \leq t \leq 4 \tau}w_{\beta}(t)}\right)<+\infty
\end{equation}

If $\tau \geq 1$, then it follows

$$
\frac{\sup _{\alpha \leq t \leq 4 \alpha} \omega_{\beta}(t)}{\inf _{\tau \leq t \leq 4 \tau} \omega_{\beta}(t)}=\frac{\sup _{\tau \leq t \leq 4 \tau} \chi(t)}{\inf _{\tau \leq t \leq 4 \tau} \chi(t)} \leq M,
$$

and consequently,

\begin{equation}\label{7.10}
    \sup _{\tau \geq 1}\left(\frac{\sup _{\tau \leq t \leq 4 \tau} \omega_{\beta}(t)}{\inf _{\tau \leq t \leq 4 \tau} \omega_{\beta}(t)}\right)<+\infty .
\end{equation}

Combining \eqref{7.8} and \eqref{7.9}, we deduce from  that there exists a constant $D_{1}>0$ independent of $x_{0}$ and $r$ such that

\begin{equation} \label{7.11}
    \frac{1}{\left|B\left(x_{0}, r\right)\right|^{\frac{N}{2}}}\left(\int_{B\left(x_{0}, r\right)}w_{\beta}(x) d x\right)\left(\int_{B\left(x_{0}, r\right)} \frac{d x}{w_{\beta}(x)}\right)^{\frac{N}{2}-1} \leq D_{1} .
\end{equation}

This finish the proof.

\section{ Some useful preliminary results }

In this section, we will derive several technical lemmas for our use later. First we begin by the radial lemma due to Lions \cite{L1}. Let $W^{1,p}(\mathbb{R}^{N})$ be the first order Sobolev space and consider the subspace of radial function namely $W^{1,p}_{rad}(\mathbb{R}^{N})$. We have
\begin{lemma}\cite{L1}\label{lemrr} Let $N\geq2$, $1\leq p<+\infty$, $u\in W^{1,p}_{rad}(\mathbb{R}^{N})$, then   there exists a positive constant $C=C(N,p)$ such that
$$|u(x)|\leq C\frac{1}{|x|^{\frac{N-1}{p}}}|u|_{p}^{\frac{p-1}{p}}|\nabla u|_{p}^{\frac{1}{p}}~~~\mbox{for}~~~~a.e x\in \mathbb{R}^{N} .$$\\
In particular, for $p=\frac{N}{2}$, we get the followin inequality:
\begin{align}\label{er1}|u(x)|\leq C\frac{1}{|x|^{\frac{2(N-1)}{N}}}|u|_{\frac{N}{2}}^{\frac{N-2}{N}}|\nabla u|_{\frac{N}{2}}^{\frac{2}{N}}~~~\mbox{for}~~~a.e~~ x\in \mathbb{R}^{N} .\end{align}
\end{lemma}
It follows that, for $N\geq 4$ , using Young inequality and the fact that $ w_{\beta}(x)\geq 1$, we get
\begin{align}\label{er2}\nonumber|u(x)|&\leq C\frac{N-2}{2}\frac{1}{|x|^{\frac{2(N-1)}{N}}}\bigg(|u|_{\frac{N}{2}}+|\nabla u|_{\frac{N}{2}}\bigg)~~~\mbox{for}~~~~a.e~~ x\in \mathbb{R}^{N}\\ \nonumber&\leq C\frac{N-2}{2}\frac{1}{|x|^{\frac{2(N-1)}{N}}}\|u\| ~~~\mbox{for}~~~~a.e~~ x\in \mathbb{R}^{N} \\&\leq C\frac{1}{|x|^{\frac{2(N-1)}{N}}}\|u\| ~~~\mbox{for}~~~~a.e~~ x\in \mathbb{R}^{N}.\end{align}
Now, we give the following Strauss compactness lemma\cite{ST}.
\begin{lemma}\label{lemst}
    Let $(P_n)_{n}$ and $(Q_n)_{n}$ be two sequences of continuous functions: $\mathbb{R}^n \rightarrow \mathbb{R}$. For $c > 0$, let $y(c) = \sup \{|t| : t = P_n(s) \text{ for some } n \text{ and } s \text{ such that } |Q_n(s)| < c |P_n(s)|\}$.

    Assume the following conditions:
    \begin{itemize}
        \item[i)] $y(c) < \infty$ for all $c > 0$. (In other words, $\displaystyle\frac{P_n}{Q_n} \rightarrow 0$ uniformly as $ n \rightarrow \infty$.)
        \item[ii)] $(u_n)_{n}$ is a sequence of measurable functions: $\mathbb{R}^N \rightarrow \mathbb{R}$ such that $\displaystyle \sup_n \int_{B} |Q_n(u_n(x))| \, dx < \infty$ for all bounded sets $B$.
        \item[iii)] $P_n(u_n(x)) \rightarrow v(x)$ for almost every $x \in \mathbb{R}^N$.
    \end{itemize}

    Then:
    \begin{itemize}
        \item[a)] $\displaystyle \int_{B} |v(x)| \, dx < \infty$ for all bounded sets $B$.
        \item[b)] Assume in addition that
        \begin{itemize}
            \item[iv)] $P_n(s) = o(Q_n(s))$ as $s \rightarrow 0$ uniformly in $n$.
            \item[v)] $u_n(x) \rightarrow 0$ as $|x| \rightarrow \infty$ uniformly in $x$ and $n$.
        \end{itemize}
        Then $\displaystyle \int_{\mathbb{R}^N} |P_n(u_n) - v| \, dx \rightarrow 0$.
    \end{itemize}
\end{lemma}
We denote by $B$ the unit ball of $\mathbb{R}^{N}$ and consider the subspace
$$ W_{0,rad}^{2,\frac{N}{2}}(B,w)=closure\{u\in
C_{0,rad}^{\infty}(B)~~|~~\displaystyle\int_{B}\big(\log(\frac{e}{|x|}\big)^{\beta(\frac{N}{2}-1)}|\Delta u|^{\frac{N}{2}}dx <\infty\}.$$  We have the following results.
\begin{lemma}\label{lem1}Let $u$ be a radially symmetric
 function in $C_{0}^{2}(\mathbb{R}^{N})$. Then, we have

\item[(i)]  Let $u$ be a radially symmetric
 function in $C_{0,rad}^{\infty}(B)$. Then, we have\begin{itemize}\item[$(i)$]\cite{ZZ}
 $$|u(x)|\leq \bigg(\displaystyle\frac{N}{\alpha_{\beta}}\big(|\log(\frac{e}{|x|}|-1\big)\bigg)^{\frac{1}{\gamma}}\displaystyle \big(\int_B w_{\beta}(x)|\Delta u|^{\frac{N}{2}}dx \big)^{\frac{2}{N}}\leq \displaystyle \bigg(\displaystyle\frac{N}{\alpha_{\beta}}\big(|\log(\frac{e}{|x|}|-1\big)\bigg)^{\frac{1}{\gamma}}\|u\|\cdot$$
\item[$(ii)$]$\displaystyle\int_{B}e^{|u(x)|^{\gamma}}dx<+\infty,~~\forall~~u\in W^{2,\frac{N}{2}}_{0,rad}(B).$

\item[(iii)]  The following embedding
is continuous $$E\hookrightarrow L^{p}(\mathbb{R}^{N})~~\mbox{for all}~~p\geq\frac{N}{2}\cdot$$
\item[(vi)]$E$ is compactly
embedded in $L^{q}(\mathbb{R}^{N})$  for all  $q\geq\frac{N}{2}\cdot$
\end{itemize}
\end{lemma}
\textit{Proof }

$(i)$ see \cite{WZ}\\

$(ii)$ From $(i)$ and using the identity $\log(\frac{e}{|x|} )-|\log(|x|)|=1~~\forall x\in B$, we get $$|u(x)|^{\gamma}\leq \displaystyle\frac{1}{\alpha_{\beta}}\bigg||\log(\frac{e}{|x|})\bigg|\displaystyle \|u\|^{\gamma}\ \leq \frac{N}{\alpha_{\beta}}\big(1+\big|\log(|x|)\big|\big) ~~\|u\|^{\gamma}. $$ Hence, using the fact that the function $r\mapsto r^{N-1}e^{\frac{\|u\|^{\gamma}(1+|\log r|)}{\alpha_{\beta}}}$ is increasing, we get
$$\displaystyle\int_{|x|<1}e^{|u|^{\gamma}}dx\leq N V_{N}\int_{0}^{1}r^{N-1}e^{\frac{N\|u\|^{\gamma}(1+|\log r|)}{\alpha_{\beta}}} dr\leq N V_{N}e^{\frac{N\|u\|^{\gamma}}{\alpha_{\beta}}} <+\infty.$$
Then $(ii)$ follows by density.\\

  $(iii)$  Since $w_{\beta}(x)\geq1$, then by Sobolev theorem, the following embedding are continuous   $$E\hookrightarrow W^{2,\frac{N}{2}}_{rad}(\mathbb{R}^{N},w_{\beta})\hookrightarrow W^{2,\frac{N}{2}}_{rad}(\mathbb{R}^{N})  \hookrightarrow L^{q}(\mathbb{R}^{N})~~\forall q\geq \frac{N}{2}.$$
We assert that the embedding
$E \rightarrow L^{q} (\mathbb{R}^{N} )$ is compact. To do this, set $Q(s) = |s|^{q}$ and $P(s) = |s|^{q+\epsilon_{0}} + |s|^{q-\epsilon_{0}} $, where $ 0 < \epsilon_{0} < q - \frac{N}{2}$.

Clearly, $ \displaystyle\frac{Q(s)}{P(s)}\rightarrow 0~~\mbox{as}~~ |s|\rightarrow+\infty$, and $\displaystyle\frac{Q(s)}{P(s)}\rightarrow 0~~\mbox{as}~~ |s|\rightarrow0$. Let $(u_{n})_{n} \in E$ be such that $u_{n}\hookrightarrow 0 $ weakly in $E$ and $u_{n}(x) \rightarrow 0~~ a.e. x \in \mathbb{R}^{N}$ . By the continuity of the embedding $E \hookrightarrow L^{q+\varepsilon_{0}} (\mathbb{R}^{N} )$ and $E \hookrightarrow L^{q-\varepsilon_0}(\mathbb{R}^N)$, we obtain that $$  \sup_n \int_{\mathbb{R}^N}|P(u_n)| < +\infty.$$ On the other hand, by  (\ref{er2}),  $u_n(x) \rightarrow 0 \text{ as } |x| \rightarrow +\infty$, uniformly in $ n \in \mathbb{N}$. Therefore, we can apply the compactness Strauss Lemma \ref{lemst}  to deduce that  $Q(u_n) \rightarrow 0 \text{ strongly in } L^1(\mathbb{R}^N)$.\\

 This concludes the lemma.
 \begin{remark}By Lemma \ref{lem1} $(ii)$ and (\ref{er2}), we have \begin{align*} \int_{\mathbb{R}^{N}}|u(x)|^{p} dx&= \int_{B}|u(x)|^{p} dx+ \int_{\mathbb{R}^{N}\setminus B}|u(x)|^{p} dx\\&\leq N V_{N}\|u\|^{p}\int_{0}^{1}r^{N-1}(1+|\log r|)^{\frac{p}{\gamma}} dr+CNV_{N}\|u\|^{p}\int_{1}^{\infty}r^{N-1-2\frac{p(N-1)}{N}}dr\\&\leq  V_{N}\|u\|^{p}+CNV_{N}\|u\|^{p}\int_{1}^{\infty}r^{N-1-2\frac{p(N-1)}{N}}dr.\end{align*}

  The last integral is finite provided $p>\frac{N^{2}}{2(N-1)}>\frac{N}{2}$. The result of the previous lemma is thus partially found.
 \end{remark}
 \begin{lemma}\cite{FMR}\label{lem2}
 Let $\Omega\subset \mathbb{R^{N}}$ be a bounded domain and $f:\overline{\Omega}\times\mathbb{R}$
  a continuous function. Let $(u_{n})_{n}$ be a sequence in $L^{1}(\Omega)$
converging to $u$ in $L^{1}(\Omega)$. Assume that $\displaystyle f(x,u_{n})$ and
$\displaystyle f(x,u)$ are also in $ L^{1}(\Omega)$. If
$$\displaystyle\int_{\Omega}|f(x,u_{n})u_{n}|dx \leq C,$$ \\where $C$
is a positive constant, then $$f(x,u_{n})\rightarrow
f(x,u)~~\mbox{in}~~L^{1}(\Omega).$$
\end{lemma}

\section{ Proof of Theorem \ref{th1.4} }

 We begin by proving the first statement of Theorem \ref{th1.4}. We have
for all $u\in E$,\begin{equation}\label{eq:1}
\int_{\R^N}(e^{\vert u\vert^{\gamma}}-1)dx = \int_{\vert x\vert\geq 1}(e^{\vert u\vert^{\gamma}}-1)dx + \int_{\vert x\vert<1}(e^{\vert u\vert^{\gamma}}-1)dx.
\end{equation}
On the one hand,
\begin{equation}\label{eq:2}
 \int_{\vert x\vert\geq 1}(e^{\vert u\vert^{\gamma}}-1)dx =\sum_{k=1}^{+\infty}\frac{1}{k!} \int_{\vert x\vert\geq 1} \vert u\vert^{\gamma k}dx.
\end{equation}
From Lemma \ref{lemrr}, we get
\begin{align}\label{eq:3}
\displaystyle \int_{\vert x\vert \geq 1} \vert u\vert^{\gamma k}dx\leq\displaystyle NV_{N}\|u\|^{\gamma k} \int_{1}^{+\infty}\frac{1}{r^{1-N+2 \gamma k\frac{N-1}{N}}}dr&\nonumber=NV_{N}\|u\|^{\gamma k}\frac{1}{-N+2 \gamma k\frac{N-1}{N}}\\ &\leq NV_{N}\|u\|^{\gamma k}\frac{1}{-N+2 \gamma \frac{N-1}{N}} ,\\ \nonumber \quad \mbox{for all}\; k\geq 1; \frac{N^{2}-4N+2}{N(N-2)}<\beta<1.
\end{align}
Combining \eqref{eq:2} and \eqref{eq:3}, we have
\begin{equation}\label{eq:4}
\int_{\vert x\vert\geq 1}(e^{\vert u\vert^{\gamma}}-1)dx \leq \frac{NV_{N}}{-N+2 \gamma \frac{N-1}{N}} \sum_{k=1}^{+\infty}\frac{\|u\|^{\gamma k}}{k!}=\frac{NV_{N}}{-N+2 \gamma \frac{N-1}{N}}  e^{\|u\|^{\gamma}}< + \infty.
\end{equation}
Now we are going to estimate the second integral in \eqref{eq:1}.  Set
\begin{equation}\label{eq:5}v(x)=\left\{\begin{array}{lll}
u(x) - u(e_1), &0\leq \vert x\vert<1,\\
0, & \vert x\vert\geq 1,
\end{array}\right.\end{equation}
where $e_1=(1,0,0,0,..........0)\in \R^N.$ Clearly $v \in W^{2,\frac{N}{2}}_{0,rad}(B,w_{\beta})$.\\
For all $\e>0$, we have
$$\vert u\vert^{\gamma}=\vert v+ u(e_1)\vert^{\gamma}\leq (1+\e)\vert v\vert^{\gamma} + \Big( 1- \frac{1}{(1+\e)^{\frac{1}{\gamma-1}}}\Big)^{1-\gamma} \vert u(e_1)\vert^{\gamma}.$$
\\
Then, from Lemma \ref{lem1} $(ii)$, we have
\begin{equation}\label{eq:6}\begin{array}{lll}
\displaystyle \int_{\vert x\vert < 1}  e^{\vert u\vert^{\gamma}}dx &\displaystyle \leq
\int_{\vert x\vert < 1}  e^{(1+\e)\vert v\vert^{\gamma}} e^{\Big( 1- \frac{1}{(1+\e)^{\frac{1}{1-\gamma}}}\Big)^{1-\gamma}\vert u(e_1)\vert^{\gamma}}dx \\
&\displaystyle \leq e^{\Big( 1- \frac{1}{(1+\e)^{\frac{1}{1-\gamma}}}\Big)^{1-\gamma}\vert u(e_1)\vert^{\gamma}}
\int_{\vert x\vert < 1}  e^{(1+\e)\vert v\vert^{\gamma}} dx<+\infty
\end{array}\end{equation}
Combining \eqref{eq:1}, \eqref{eq:4}, \eqref{eq:6}, (\ref{er2}) and Lemma \ref{lem1} $(ii)$, we conclude that
$$\int_{\R^N}(e^{\vert u\vert^{\gamma}}-1)dx< +\infty, \quad \mbox{for all}\, u \in E.$$
This ends the proof of the first item .\\ By \eqref{eq:4} we have
\begin{align}\label{eq:4.7}\sup_{u\in E, \Vert u\Vert \leq 1 }\int_{|x|\geq1}(e^{\alpha \vert u\vert^{\gamma}}-1)dx\leq \sup_{u\in E, \Vert u\Vert \leq 1}\frac{NV_{N}}{-N+2 \gamma \frac{N-1}{N}}  e^{\|u\|^{\gamma}} \leq \frac{NV_{N}}{-N+2 \gamma \frac{N-1}{N}}  e .\end{align}
On the other hand, by \eqref{eq:6}, (\ref{er2}) and using the radial lemma \ref{lem1}$(i)$, we get \begin{align}\label{eq:8}\sup_{u\in E, \Vert u\Vert \leq 1 }\int_{|x|\leq1}(e^{\alpha \vert u\vert^{\gamma}}-1)dx &\nonumber\leq e^{\Big( 1- \frac{1}{(1+\e)^{\frac{1}{1-\gamma}}}\Big)^{1-\gamma}\vert u(e_1)\vert^{\gamma}}\sup_{u\in E, \Vert u\Vert \leq 1}\int_{\vert x\vert < 1}  e^{(1+\e)\vert v\vert^{\gamma}} dx\\& \nonumber\leq e^{\Big( 1- \frac{1}{(1+\e)^{\frac{1}{1-\gamma}}}\Big)^{1-\gamma}(C\|u\|))^{\gamma}}\sup_{u\in E, \Vert u\Vert \leq 1}\int_{\vert x\vert < 1}  e^{(1+\e)\vert v\vert^{\gamma}} dx\\&\leq e^{\Big( 1- \frac{1}{(1+\e)^{\frac{1}{1-\gamma}}}\Big)^{1-\gamma}(C)^{\gamma}}\sup_{u\in E, \Vert u\Vert \leq 1}\int_{\vert x\vert < 1}  e^{(1+\e)\vert v\vert^{\gamma}} dx.\end{align} Let $\alpha < \alpha_\beta$. Clearly, there exists $\e>0$ such that $\alpha(1+\e)<\alpha_\beta$.\\ Do not forgot that
\begin{align}\label{eq:7}
\Vert  v\Vert^\frac{N}{2}_{W^{2,\frac{N}{2}}_{0,rad}(B)}&=\int_B \vert \Delta v\vert^\frac{N}{2}\left(\log \left(\frac{e}{|x|}\right)\right)^{\beta(\frac{N}{2}-1)}  dx\\ \nonumber&= \int_B \vert \Delta u\vert^\frac{N}{2}w_{\beta}(x) dx\leq \Vert u\Vert^\frac{N}{2}\leq 1.
\end{align}
Then, $$\sup_{u\in E, \Vert u\Vert\leq 1 }
\int_{\vert x\vert < 1}  e^{\alpha(1+\e)\vert v\vert^{\gamma}} dx\leq \sup\big\{\int_{\vert x\vert < 1}  e^{\alpha(1+\e)\vert v\vert^{\gamma}} dx,~~v\in W^{2,\frac{N}{2}}_{0,rad}(w,B),~~\|v\|_{W^{2,\frac{N}{2}}_{0,rad}(B)}\leq1\big\}.$$

So by (\ref{eq:8}), there exists $C(\beta)>0$ such that

\begin{equation}\label{eq:10}
\displaystyle \sup_{u\in E, \Vert u\Vert  \leq 1 }\int_{\vert x\vert < 1}  e^{\alpha\vert u\vert^{\gamma}}dx
\displaystyle \leq e^{\Big( 1- \frac{1}{(1+\e)^{\frac{1}{1-\gamma}}}\Big)^{1-\gamma}(C)^{\gamma}}C(\beta).
\end{equation}
Combining \eqref{eq:8} and \eqref{eq:7}, we get
$$\sup_{u\in E, \Vert u\Vert \leq 1 }\int_{|x|<1}(e^{\alpha \vert u\vert^{\gamma}}-1)dx<+\infty.$$
Furthermore
\begin{equation}\label{eq:9}
\int_{\vert x\vert\geq 1}(e^{\alpha\vert u\vert^{\gamma}}-1)dx =\sum_{k=1}^{+\infty}\frac{\alpha^k}{k!} \int_{\vert x\vert\geq 1} \vert u\vert^{\gamma k}dx.
\end{equation}
Combining \eqref{eq:3} and \eqref{eq:9}, we infer
\begin{equation}\label{eq:10}\sup_{u\in E, \Vert u\Vert\leq 1 }\int_{\vert x\vert\geq 1}(e^{\alpha \vert u\vert^{\gamma}}-1)dx< +\infty.\end{equation}
It follows from \eqref{eq:8} and \eqref{eq:10} that
$$\sup_{u\in E, \Vert u\Vert \leq 1 }\int_{\R^N}(e^{\alpha \vert u\vert^{\gamma}}-1)dx< +\infty, \quad \mbox{for all}\, \alpha <\alpha_\beta.$$
In the next step , we show that if $\alpha>\alpha_{\beta}$, then the
supremum is infinite.
Now, we will use particular functions \cite{WZ}, namely the Adams' functions.
  We consider the sequence defined for all $n\geq3$  by
 \begin{equation}\label{eq:5.2}w_{n}(x)=C(N,\beta)\displaystyle \left\{
      \begin{array}{rclll}
&\displaystyle\bigg(\frac{\log (e\sqrt[N]{n} )}{\alpha_{\beta}}\bigg)^{\frac{1}{\gamma}}-\frac{|x|^{2(1-\beta)}}{2\big(\alpha_{\beta}\big)^{\frac{1}{\gamma}}
\big(\log  (e\sqrt[N]{n} )\big)^{\frac{\gamma-1}{\gamma}}}\\&+\frac{1}{2\big(\alpha_{\beta}\big)^{\frac{1}{\gamma}}(\frac{1}{n})^{\frac{2(1-\beta)}{N}}
\big(\log  (e\sqrt[N]{n} )\big)^{\frac{\gamma-1}{\gamma}}}& \mbox{ if } 0\leq |x|\leq \frac{1}{\sqrt[N]{n}}\\\\
       &\displaystyle  \frac{1}{\alpha^{\frac{1}{\gamma}}_{\beta}\big(\log e\sqrt[N]{n} \big)^{\frac{2(1-\beta)}{N}}}\bigg(\log(\frac{e}{|x|}\bigg)^{1-\beta} & \mbox{ if } \frac{1}{\sqrt[N]{n}}\leq|x|\leq \frac{1}{2}\\
       &\zeta_{n}&\mbox{ if}~~ |x| \geq \frac{1}{2}
 \end{array}
    \right.
  \end{equation}
     where $\displaystyle C(N,\beta)=\frac{(\frac{1}{2})^{\frac{2}{N}}\alpha^{\frac{1}{\gamma}}_{\beta}}{V^{\frac{2}{N}}_{N}(1-\beta)^{1-\frac{2}{N}}(N-2)} $, $\zeta_{n}\in C^{\infty}_{0,rad}(B)$ is such that\\
      $\displaystyle\zeta_{n}(\frac{1}{2})= \frac{1}{\alpha^{\frac{1}{\gamma}}_{\beta}\big(\log e\sqrt[N]{n} \big)^{\frac{2(1-\beta)}{N}}}\big(\log 2e \big)^{1-\beta}$,
      $\displaystyle\frac{\partial \zeta_{n}}{\partial r}(\frac{1}{2}) = \frac{-2(1-\beta)}{\alpha^{\frac{1}{\gamma}}_{\beta}\big(\log e\sqrt[N]{n} \big)^{\frac{2(1-\beta)}{N}}} \big(\log (2e)\big)^{-\beta}$ \newline  $\displaystyle\zeta_{n}(1)=\frac{\partial \zeta_n}{\partial r}(1)=0$ and $\xi_{n}$, $\nabla \xi_{n}$, $\Delta \xi_{n}$ are all $\displaystyle o\bigg(\frac{1}{[\log (e\sqrt[N]{n})]^{\frac{1}{\gamma}}}\bigg)$. Here,  $\displaystyle\frac{\partial \zeta_{n}}{\partial r}$ denotes the first derivative of $\zeta_{n}$ in the radial variable $r=|x|$. \\
     Let $v_{n}(x)=\displaystyle\frac{w_{n}}{\|w_{n}\|}$.
 We have, $v_{n}\in E$ , $\|v_{n}\|^{\frac{N}{2}}=1.$ \\

 We compute $\Delta w_{n}(x)$, we get \begin{equation*}\label{eq:4.2}\Delta w_{n}(x)=C(N,\beta)\displaystyle \left\{
      \begin{array}{rclll}
\displaystyle\frac{-2(1-\beta)(N-2\beta)|x|^{-2\beta}}{\alpha^{\frac{1}{\gamma}}_{\beta}\big(\log e\sqrt[N]{n} \big)^{\frac{2(1-\beta)}{N}}\big(\log  (e\sqrt[N]{n})\big)^{\frac{\gamma-1}{\gamma}}} & \mbox{ if } 0\leq |x|\leq \frac{1}{\sqrt[N]{n}}\\\\
       \displaystyle \frac{-(1-\beta)}{|x|^{2}} \frac{\bigg(\log(\frac{e}{|x|})\bigg)^{-\beta}\bigg((N-2)+\beta\big(\log \frac{e}{|x|}\big)^{-1}\bigg)}{\alpha^{\frac{1}{\gamma}}_{\beta}\big(\log e\sqrt[N]{n} \big)^{\frac{2(1-\beta)}{N}}} & \mbox{ if } \frac{1}{\sqrt[4]{n}}\leq|x|\leq \frac{1}{2}\\
       \Delta \zeta_{n} &\mbox{ if} |x| \geq \frac{1}{2}
 \end{array}
    \right.
  \end{equation*}
  So, \begin{align*}    \frac{1}{C^{\frac{N}{2}}(N,\beta)} \|\Delta w_{n}\|_{\frac{N}{2},w}^{\frac{N}{2}}    = &\underbrace{N V_{N}\int^{\frac{1}{    \sqrt[N]{n}}}_{0}r^{N-1}|\Delta w_{n}(x)|^{\frac{N}{2}}\big(\log \frac{e}{r}\big)^{\beta(\frac{N}{2}-1)}dr}_{I_{1}}  \\
  & +\underbrace{N V_{N}\int^{\frac{1}{2}}_{\frac{1}{\sqrt[N]{n}}}r^{N-1}|\Delta w_{n}(x)|^{\frac{N}{2}}\big(\log \frac{e}{r}\big)^{\beta(\frac{N}{2}-1)}dr}_{I_{2}}    \\  &+\underbrace{N V_{N}\int^{1}_{\frac{1}{2}}r^{N-1}|\Delta \zeta_{n}(x)|^{\frac{N}{2}}\big(\log \frac{e}{r}\big)^{\beta(\frac{N}{2}-1)}dr+ N V_{N}\int^{+\infty}_{1}|\triangle \zeta_{n}|^{\frac{N}{2}}\chi(r)r^{N-1}dr}_{I_{3}}\end{align*}
By using integration by parts, we obtain,$$\begin{array}{rclll}\displaystyle I_{1}&= & \displaystyle NV_{N}\frac{2^{\frac{N}{2}}(1-\beta)^{\frac{N}{2}}(N-2\beta)^{\frac{N}{2}}}{\big(\alpha_{\beta}\big)^{\frac{N}{2\gamma}}(\log  (e\sqrt[N]{n})\big)^{1-\beta}\big(\log  (e\sqrt[N]{n})\big)^{\frac{N(\gamma-1)}{2\gamma}}}\int^{\frac{1}{\sqrt[N]{n}}}_{0} r^{N(1-\beta)-1}\big(\log \frac{e}{r}\big)^{\beta(\frac{N}{2}-1)}dr\\
&=&\displaystyle NV_{N}\frac{2^{\frac{N}{2}}(1-\beta)^{\frac{N}{2}}(N-2\beta)^{\frac{N}{2}}}{\big(\alpha_{\beta}\big)^{\frac{N}{2\gamma}}(\log  (e\sqrt[N]{n})\big)^{1-\beta}\big(\log  (e\sqrt[N]{n})\big)^{\frac{N(\gamma-1)}{2\gamma}}}\left[  \frac{r^{N(1-\beta)}}{N(1-\beta)}(\log \frac{e}{r}\big)^{\beta(\frac{N}{2}-1)} \right]^{\frac{1}{\sqrt[N]{n}}}_{0}\\&+& NV_{N}\frac{\beta(\frac{N}{2}-1)2^{\frac{N}{2}}(1-\beta)^{\frac{N}{2}}(N-2\beta)^{\frac{N}{2}}}{\big(\alpha_{\beta}\big)^{\frac{N}{2\gamma}}(\log  (e\sqrt[N]{n})\big)^{1-\beta}\big(\log  (e\sqrt[N]{n})\big)^{\frac{N(\gamma-1)}{2\gamma}}}\displaystyle\int^{\frac{1}{\sqrt[N]{n}}}_{0}\frac{r^{N(1-\beta)}-1}{N(1-\beta)} \big(\log \frac{e}{r}\big)^{\beta(\frac{N}{2}-1)-1} dr\\
&=&\displaystyle o\big(\frac{1}{\log e\sqrt[N]{n}}\big)\cdot\\

\end{array}$$
Also,\begin{align}\nonumber\displaystyle \frac{2}{N}C^{\frac{N}{2}}(N,\beta)I_{2}=\\\nonumber= &\displaystyle C^{\frac{N}{2}}(N,\beta) NV_{N}\frac{(1-\beta)^{\frac{N}{2}}}{\big(\alpha_{\beta}\big)^{\frac{N}{2\gamma}}\big(\log  (e\sqrt[N]{n})\big)^{1-\beta}}\int_{\frac{1}{\sqrt[N]{n}}}^{\frac{1}{\frac{1}{2}} }\frac{1}{r}\big(\log \frac{e}{r}\big)^{-\beta}\big ((N-2)+\beta \big( \log \frac{e}{r}\big)^{-1}\big)^{\frac{N}{2}}dr\\\nonumber
=& \nonumber\displaystyle C^{\frac{N}{2}}(N,\beta)NV_{N}\frac{(1-\beta)^{\frac{N}{2}}(N-2)^{\frac{N}{2}}}{\big(\alpha_{\beta}\big)^{\frac{N}{2\gamma}}\big(\log  (e\sqrt[N]{n})\big)^{1-\beta}}\int_{\frac{1}{\sqrt[N]{n}}}^{\frac{1}{\frac{1}{2}} }\frac{1}{r}\big(\log \frac{e}{r}\big)^{-\beta}\big (1+o( \log \frac{e}{r}\big)^{-1}\big)^{\frac{N}{2}}dr\\=& \nonumber\displaystyle C^{\frac{N}{2}}(N,\beta) NV_{N}\frac{(1-\beta)^{\frac{N}{2}}(N-2)^{\frac{N}{2}}}{\big(\alpha_{\beta}\big)^{\frac{N}{2\gamma}}\big(\log  (e\sqrt[N]{n})\big)^{1-\beta}}\bigg(\int_{\frac{1}{\sqrt[N]{n}}}^{\frac{1}{\frac{1}{2}} }\frac{1}{r}\big(\log \frac{e}{r}\big)^{-\beta}dr +\int_{\frac{1}{\sqrt[N]{n}}}^{\frac{1}{\frac{1}{2}} }\frac{1}{r}\big(\log \frac{e}{r}\big)^{-\beta}o( \log \frac{e}{r}\big)^{-1}dr\bigg)\\=& \nonumber\displaystyle C^{\frac{N}{2}}(N,\beta)NV_{N}\frac{(1-\beta)^{\frac{N}{2}}(N-2)^{\frac{N}{2}}}{\big(\alpha_{\beta}\big)^{\frac{N}{2\gamma}}\big(\log  (e\sqrt[N]{n})\big)^{1-\beta}}\left [ \frac{1}{1-\beta}\big( \log \frac{e}{r}\big)^{1-\beta}\right]^{\frac{1}{\sqrt[N]{n}}}_{\frac{1}{2}}\\ \nonumber&-
C^{\frac{N}{2}}(N,\beta)NV_{N}\frac{(1-\beta)^{\frac{N}{2}}(N-2)^{\frac{N}{2}}}{\big(\alpha_{\beta}\big)^{\frac{N}{2\gamma}}\big(\log  (e\sqrt[N]{n})\big)^{1-\beta}}\bigg(\int_{\frac{1}{\sqrt[N]{n}}}^{\frac{1}{\frac{1}{2}} } \frac{1}{r}\big(\log \frac{e}{r}\big)^{-\beta}o( \log \frac{e}{r}\big)^{-1}dr\bigg)\\\nonumber
=& \nonumber\displaystyle 1+ o\big(\frac{1}{(\log e\sqrt[N]{n})^{1-\beta}}\big)\cdot
\end{align}
and $I_{3}=\displaystyle  o\big(\frac{1}{(\log e\sqrt[4]{n})^{\frac{2}{\gamma}}}\big).$ Then $\frac{2}{N}\|\Delta w_{n}\|_{\frac{N}{2},w}^{\frac{N}{2}}=1+o\big(\frac{1}{(\log e\sqrt[N]{n})^{\frac{N}{\gamma}}}\big)$.\\
In the sequel we prove the following key lemma.
\begin{lemma}\label{lem6} The Adams' function given by (\ref{eq:5.2}) verifies
    $\displaystyle \lim_{n\rightarrow+\infty}\|w_{n}\|^{\frac{N}{2}}=1.$

\end{lemma}
    {\it Proof}~~We have
   \begin{align}\nonumber\|w_{n}\|^{\frac{N}{2}}&=  \displaystyle\frac{2}{N}\int_{B} w_{\beta}(x)|\Delta w_{n}|^{\frac{N}{2}}dx+\frac{2}{N}\int_{B}|\nabla w_{n}|^{\frac{N}{2}}dx+\frac{2}{N}\int_{B}~~|w_{n}|^{\frac{N}{2}}dx \\\nonumber
&=\displaystyle 1+o\big(\frac{1}{(\log e\sqrt[N]{n})^{\frac{2}{\gamma}}}\big)+\nonumber\frac{2}{N}\int_{0\leq |x|\leq \frac{1}{\sqrt[N]{n}}}|w|^{\frac{N}{2}}_{n}dx+\frac{2}{N}\int_{\sqrt[N]{n}\leq |x|\leq \frac{1}{2}} |w_{n}|^{\frac{N}{2}}dx\\ \nonumber &+\displaystyle\frac{2}{N}\int_{|x|\geq \frac{1}{2}}|\zeta|^{\frac{N}{2}}_{n}dx\\&\nonumber+\displaystyle\underbrace{\frac{2}{N}\int_{0\leq |x|\leq \frac{1}{\sqrt[N]{n}}}| \nabla w_{n}|^{\frac{N}{2}}dx}_{I_{1}'}+\underbrace{\frac{2}{N}\int_{\sqrt[N]{n}\leq |x|\leq \frac{1}{2}}| \nabla w_{n}|^{\frac{N}{2}}dx}_{I_{2}'}+\underbrace{\frac{2}{N}\int_{|x|\geq \frac{1}{2}}|\nabla \zeta_{n}|^{\frac{N}{2}}dx}_{I_{3}'}\cdot
\end{align}
We have,$$\begin{array}{rclll}\displaystyle I'_{1}&= & \displaystyle 2V_{N}\frac{C^{\frac{N}{2}}(N,\beta)(1-\beta)^{\frac{N}{2}}}{\alpha^{\frac{N}{2\gamma}}_{\beta}\big(\log  (\sqrt[N]{n})\big)^{\frac{N(\gamma-1)}{2\gamma}}}\int^{\frac{1}{\sqrt[N]{n}}}_{0} r^{N(2-\beta)-1}dr\\
&=&\displaystyle 2V_{N}\frac{C^{\frac{N}{2}}(N,\beta)(1-\beta)^{\frac{N}{2}}}{\alpha^{\frac{N}{2\gamma}}_{\beta}\big(\log  (\sqrt[N]{n})\big)^{\frac{N(\gamma-1)}{2\gamma}}}\left[  \frac{r^{N(2-\beta)}}{N(2-\beta)} \right]^{\frac{1}{\sqrt[N]{n}}}_{0}
\\&=&\displaystyle 2V_{N}\frac{C^{\frac{N}{2}}(N,\beta)(1-\beta)^{\frac{N}{2}}}{\alpha^{\frac{N}{2\gamma}}_{\beta}n^{2-\beta}N(2-\beta)\log  (\sqrt[N]{n})\big)^{\frac{N(\gamma-1)}{2\gamma}}}\\
&=&\displaystyle o\big(\frac{1}{n^{2-\beta}\log e\sqrt[N]{n}}\big)\cdot

\end{array}$$
Also, using the fact that the function $r\mapsto r^{\frac{N}{2}-1}\big(\log \frac{e}{r}\big)^{-\frac{N}{2}\beta}$ is increasing on $[0,1]$, we get  $$\begin{array}{rclll}\displaystyle I'_{
2}&= & \displaystyle 2V_{N}C^{\frac{N}{2}}(N,\beta)\frac{(1-\beta)^{\frac{N}{2}}}{\big(\alpha_{\beta}\big)^{\frac{N}{2\gamma}}\big(\log  (e\sqrt[N]{n})\big)^{1-\beta}}\int_{\frac{1}{\sqrt[N]{n}}}^{\frac{1}{2}} r^{\frac{N}{2}-1}\big(\log \frac{e}{r}\big)^{-\frac{N}{2}\beta}dr\\
&\leq & \displaystyle 2V_{N}C^{\frac{N}{2}}(N,\beta)\frac{(1-\beta)^{\frac{N}{2}}}{\big(\alpha_{\beta}\big)^{\frac{N}{2\gamma}}\big(\log  (e\sqrt[N]{n})\big)^{1-\beta}}(\frac{1}{2})^{\frac{N}{2}-1}\big(\log 2e\big)^{-\frac{N}{2}\beta}\\
&=&\displaystyle o\bigg(\frac{1}{[\log (e\sqrt[N]{n})]^{1-\beta}}\bigg)\\
\end{array}$$
andand $I'_{3}=\displaystyle  o\big(\frac{1}{(\log e\sqrt[N]{n})^{\frac{2}{\gamma}}}\big).$
    For $\displaystyle|x|\leq \frac{1}{\sqrt[N]{n}}$, $$\displaystyle |w_{n}|^{\frac{N}{2}}\leq C(N,\beta)^{\frac{N}{2}} \bigg(\frac{\log (e\sqrt[N]{n} )}{\alpha_{\beta}}\bigg)^{\frac{1}{\gamma}}+ \frac{1}{2\big(\alpha_{\beta}\big)^{\frac{1}{\gamma}}(\frac{1}{n})^{\frac{2(1-\beta)}{N}}
\big(\log  (e\sqrt[N]{n} )\big)^{\frac{\gamma-1}{\gamma}}}\bigg)^{\frac{N}{2}}\cdot$$Then,
  \begin{align}\nonumber&\int_{0\leq |x|\leq \frac{1}{\sqrt[N]{n}}}|w_{n}|^{\frac{N}{2}}dx\leq\\\nonumber & NV_{N} C(N,\beta)^{\frac{N}{2}}\bigg(\bigg(\frac{\log (e\sqrt[N]{n} )}{\alpha_{\beta}}\bigg)^{\frac{1}{\gamma}}+ \frac{1}{2\big(\alpha_{\beta}\big)^{\frac{1}{\gamma}}(\frac{1}{n})^{\frac{2(1-\beta)}{N}}
\big(\log  (e\sqrt[N]{n} )\big)^{\frac{\gamma-1}{\gamma}}}\bigg)^{\frac{N}{2}}\int^{\frac{1}{\sqrt[N]{n}}}_{0} r^{N-1}dr=o_{n}(1)\end{align}
 Also, Also,$$\int_{\frac{1}{\sqrt[N]{n}}\leq |x|\leq \frac{1}{2}}|w_{n}|^{\frac{N}{2}}dx= \frac{NV_{N} C(N,\beta)^{\frac{N}{2}} }{\big(\alpha_{\beta}\big)^{\frac{N}{2\gamma}}\big(\log  (e\sqrt[N]{n})\big)^{1-\beta}}\int^{\frac{1}{2}}_{\frac{1}{\sqrt[N]{n}}}r^{N-1}\big(\log(\frac{e}{r}\big)\big)^{\frac{2(1-\beta)}{N}}dr$$
 Using the fact that the function $r\mapsto r^{N-1}\big(\log \frac{e}{r}\big)^{\frac{2(1-\beta)}{N}}$ is increasing on $[0,1]$,   we obtain $$\int_{\frac{1}{\sqrt[N]{n}}\leq |x|\leq \frac{1}{2}}|w_{n}|^{\frac{N}{2}}dx\leq \frac{NV_{N} C(N,\beta)^{\frac{N}{2}} }{\big(\alpha_{\beta}\big)^{\frac{N}{2\gamma}}\big(\log  (e\sqrt[N]{n})\big)^{1-\beta}}\frac{1}{2^{N-1}}\big(\log  (2e)\big)^{\frac{2(1-\beta)}{N}}=o_{n}(1).$$
 Finaly,$$\int_{|x|\geq {\frac{1}{2}}}  |w_{n}|^{\frac{N}{2}}dx=\int_{|x|\geq \frac{1}{2}}|\zeta_{n}|^{\frac{N}{2}}dx= o_{n}(1)$$
 Then, $\|\ w_{n}\|^{\frac{N}{2}}=1+o\big(\frac{1}{(\log e\sqrt[N]{n})^{\frac{2}{\gamma}}}\big)$.
The Lemma is proved. \\ Now, let $v_{n}(x)=\displaystyle\frac{w_{n}}{\|w_{n}\|}$. From the definition of $w_{n}$, it is easy to see that  $$-\frac{|x|^{2(1-\beta)}}{2\big(\alpha_{\beta}\big)^{\frac{1}{\gamma}}
\big(\log  (e\sqrt[N]{n} )\big)^{\frac{\gamma-1}{\gamma}}}+\frac{1}{2\big(\alpha_{\beta}\big)^{\frac{1}{\gamma}}(\frac{1}{n})^{\frac{2(1-\beta)}{N}}
\big(\log  (e\sqrt[N]{n} )\big)^{\frac{\gamma-1}{\gamma}}}\geq0 ~~\mbox{for all}~~ 0\leq |x|\leq \frac{1}{\sqrt[N]{n}}\cdot$$ Then, for all $0\leq |x|\leq \frac{1}{\sqrt[N]{n}}$, $\displaystyle| w_{n}|^{\frac{N}{2}}\geq \displaystyle\displaystyle\bigg(\frac{\log (e\sqrt[N]{n} )}{\alpha_{\beta}}\bigg)^{\frac{N}{2\gamma}}\cdot$ Let $\overline{\alpha}=\frac{\alpha}{\alpha_{\beta}},$  we have
$$\begin{array}{llll}\displaystyle
\sup_{u\in E, \Vert u\Vert \leq 1 }\int_{\R^N}(e^{\alpha\vert u\vert^{\gamma}}-1)dx &\displaystyle\geq \lim_{n\longrightarrow +\infty }\int_{\vert x\vert \leq \frac{1}{\sqrt[N]{n}} }(e^{\alpha |v_{n}|^{\gamma}}-1)dx\\ &\displaystyle\geq \lim\limits_{n\longrightarrow +\infty }NV_{N}\int^{\frac{1}{\sqrt[N]{n}}}_{0}\big(r^{N-1}e^{\overline{\alpha}\log\big(e\sqrt[N]{n}\big)}-r^{N-1}\big)dr\\&\displaystyle\geq \lim\limits_{n\longrightarrow \infty } V_{N}\frac{1}{n}\big(e^{\sqrt[N]{n}}-1 \big)=+\infty~~\mbox{if}~~\overline{\alpha}>1.
 \end{array}$$Then,
 $$\sup_{u\in E, \Vert u\Vert \leq 1 }\int_{\R^N}(e^{\alpha \vert u\vert^{\gamma}}-1)dx= +\infty~~\forall~~\alpha\geq\alpha_{\beta}.$$

\section{A Lions-type compactness concentration lemma}
In the sequel, we prove a concentration compactness result of Lions type.
\begin{lemma}\label{Lionstype} Let
$(u_{k})_{k}$  be a sequence in $E$. Suppose that, \newline $\|u_{k}\|=1$, $u_{k}\rightharpoonup u$ weakly in $E$, $u_{k}(x)\rightarrow u(x) ~~a.e~x\in \mathbb{R}^{N}$, $\nabla u_{k}(x)\rightarrow\nabla u(x) ~~a.e~x\in \mathbb{R}^{N}$, $\Delta u_{k}(x)\rightarrow\Delta u(x) ~~a.e~x\in \mathbb{R}^{N}$ and $u\not\equiv 0$. Then
$$\displaystyle\sup_{k}\int_{B}(e^{p~\alpha_{\beta}
|u_{k}|^{\gamma}}-1)dx< +\infty,~~\mbox{where}~~ \alpha_{\beta}=N[(N-2)N V_{N}]^{\frac{2}{(N-2)(1-\beta)}}(1-\beta)^{\frac{1}{(1-\beta)}},$$
for all $1<p<U(u)$ where $U(u)$ is given by:
 $$U(u):=\displaystyle \left\{
      \begin{array}{rcll}
&\displaystyle\frac{1}{(1-\|u\|^{\frac{N}{2}})^{\frac{2\gamma}{N}}}& \mbox{ if }\|u\| <1\\
       &+\infty& \mbox{ if } \|u\|=1\\
 \end{array}
    \right.$$
\end{lemma}
\textit{Proof}
For $a,~b \in \R,~~q>1$. If $q'$ its conjugate i.e. $\frac{1}{q}+\frac{1}{q'}=1$
we have, by Young inequality, that
$$(e^{a+b}-1)\leq \frac{1}{q}(e^{qa}-1)+ \frac{1}{q'}(e^{q'b}-1).$$

Also, we have
\begin{equation}\label {eq:3.1}
(1+a)^{q}\leq (1+\varepsilon) a^{q}+(1-\frac{1}{(1+\varepsilon)^{\frac{1}{q-1}}})^{1-q},~~\forall a\geq0,~~\forall\varepsilon>0~~\forall q>1.
\end{equation}
So, we get
$$
      \begin{array}{rcll}
|u_{k}|^{\gamma}&=& |u_{k}-u+u|^{\gamma}\\
&\leq& (|u_{k}-u|+|u|)^{\gamma}\\
&\leq& (1+\varepsilon)|u_{k}-u|^{\gamma}+\big(1-\frac{1}{(1+\varepsilon)^{\frac{1}{\gamma-1}}}\big)^{1-\gamma}|u|^{\gamma}\\
 \end{array}
   $$
which implies that
  \begin{align*}
\int_{\mathbb{R}^{4}}\big(e^{p~\alpha_{\beta}
|u_{k}|^{\gamma}}-1\big)dx &\leq
\frac{1}{q}\int_{\mathbb{R}^{N}}\big( e^{pq~\alpha_{\beta}
(1+\varepsilon)|u_{k}-u|^{\gamma}}-1\big)dx\\
&+\displaystyle\frac{1}{q'}\int_{\mathbb{R}^{N}} \bigg(e^{pq'~\alpha_{\beta}
(1-\frac{1}{(1+\varepsilon)^{\frac{1}{\gamma-1}}})^{1-\gamma}|u|^{\gamma}}-1\bigg )dx,
\end{align*}
for any $p>1$.
From Lemma \ref{lem1} $(ii)$, the last integral is finite.\\ To finish the
proof, we need to prove that for all $p$ such that $1<p<U(u)$,
\begin{equation}\label {eq:3.2}
\sup_{k}\int_{\mathbb{R}^{N}}\bigg( e^{pq~\alpha_{\beta}
(1+\varepsilon)|u_{k}-u|^{\gamma}}-1\bigg )dx<+\infty,
\end{equation}
 for some $\varepsilon>0$ and $q>1$.\\
In what follows, we assume that
$\|u\|<1$ and in the case that $\|u\|=1$, the proof  is similar.\\
When $$\|u\|<1$$
and
$$p<\displaystyle\frac{1}{(1-\|u\|^{\frac{N}{2}})^{\frac{2\gamma}{N}}},$$
there exists $\nu>0$ such that
$$p(1-\|u\|^{\frac{N}{2}})^{\frac{2\gamma}{N}}(1+\nu)<1.$$
On the other hand,
 by Brezis-Lieb's Lemma \cite{Br} we have
 \begin{equation*}\label {eq:2.4}
 \|u_{k}-u\|^{\frac{N}{2}}=\|u_{k}\|^{\frac{N}{2}}-\|u\|^{\frac{N}{2}}+o(1)~~\mbox{where}~~
  o(1)\rightarrow 0~~ \mbox{as}~~k\rightarrow +\infty.
\end{equation*}
   Then,
   \begin{align}\nonumber\|u_{k}-u\|^{\frac{N}{2}}=1-\|u\|^{\frac{N}{2}}+o(1),\end{align}
and so

$$\displaystyle\lim_{k\rightarrow+\infty}\|u_{k}-u\|^{\gamma}=(1-\|u\|^{\frac{N}{2}})^{\frac{2\gamma}{N}}.$$

 Therefore, for every
$\varepsilon>0$, there exists $k_{\varepsilon}\geq 1$ such that
$$\|u_{k}-u\| ^{\gamma}\leq (1+\varepsilon)(1-\|u\|^{\frac{N}{2}})^{\frac{2\gamma}{N}},~~\forall ~~k\geq k_{\varepsilon}.$$
If we take $q=1+\varepsilon$ with
$\varepsilon=\sqrt[3]{1+\nu}-1$, then  $\forall k\geq k_{\varepsilon}$, we have
 $$pq(1+\varepsilon)\|u_{k}-u\|^{\gamma}\leq 1.$$
Consequently,
$$\begin{array}{rlll}
\displaystyle\int_{\mathbb{R}^{N}}\bigg( e^{pq~\alpha_{\beta}
(1+\varepsilon)|u_{k}-u|^{\gamma}}-1\bigg)dx&\leq&
\displaystyle\int_{\mathbb{R}^{N}} \bigg(e^{
(1+\varepsilon)pq~\alpha_{\beta}(\frac{|u_{k}-u|}{\|u_{k}-u\|})^{\gamma}\|u_{k}-u\|^{\gamma}}-1\bigg)dx\\
&\leq&\displaystyle\int_{\mathbb{R}^{N}} \bigg(e^{~\alpha_{\beta}(\frac{|u_{k}-u|}{\|u_{k}-u\|})^{\gamma}}-1\bigg)dx\\
&\leq &\displaystyle\sup _{\|u\|\leq 1}\displaystyle\int_{\mathbb{R}^{N}}
\bigg(e^{~\alpha_{\beta}|u|^{\gamma}}-1\bigg)dx <+\infty.
\end{array}
$$
Now, (\ref{eq:3.2})  follows from  (\ref{eq:1.3}).
 This complete the proof and lemma \ref{Lionstype} is proved.
\section{The variational formulation for the problem (\ref{eq:1.1})}
Note that, by the hypothesis ($H_{3}$), for any $\varepsilon>0$, there exists $\delta_{0}>0$
   such that \begin{equation}\label{e1}|f(t)|\leq \varepsilon |t|^{\frac{N}{2}-1},~~\forall ~~0<|t|\leq \delta_{0}.
   \end{equation}
   Moreover, since $f$ is critical at infinity, for every $\varepsilon>0$, there exists $C_{\varepsilon}>0 $ such that
   \begin{equation}\label{e2}~~\forall t\geq C_{\varepsilon}~~|f(t)|\leq \varepsilon \exp( ~a|t|^{\gamma}-1)~~\mbox{with}~~a>\alpha_{0}.
   \end{equation}In particular, we obtain for $q\geq2$,\begin{equation}\label{e3} ~|f(t)|\leq \frac{\varepsilon}{C^{q-1}_{\varepsilon}}|t|^{q-1} \exp(a ~|t|^{\gamma}-1)~~\mbox{with}~~a>\alpha_{0}.
   \end{equation}
 Hence, using   (\ref{e1}), (\ref{e2}), (\ref{e3})  and the continuity of $f$,  for every $\varepsilon>0$, for every $q>N$, there exists  a positive constant $C$  such that \begin{align}\label{imp} |f(t)|\leq \varepsilon |t|^{\frac{N}{2}-1} +C |t|^{q-1}\big(e^{a ~|t|^{\gamma}}-1\big), ~~~~~~\forall\  t\in\mathbb{R}, ~~\forall~a>\alpha_{0}.\end{align}
 It follows from (\ref{imp}) and $(H_{2})$, that for all $\varepsilon>0$, there exists $C>0$ such that
\begin{equation}\label {eq:1.10}
F(t)\leq \varepsilon |t|^{\frac{N}{2}}+C |t|^{q}\big(e^{a~|t|^{\gamma}}-1\big),~~~~~~\mbox{for all}~~t,\forall~a>\alpha_{0}
\end{equation}

 So, by (\ref{eq:1.3}) and  (\ref{eq:1.10}) the functional $\mathcal{J}$ given by (\ref{energy}), is well defined. Moreover, by standard arguments,  $\mathcal{J}\in  C^{1}(E,\mathbb{R})$.

\subsection{The mountain pass geometry of the energy}

In the sequel, we prove that the functional $\mathcal{J}$ has a mountain pass
geometry.

\begin{proposition}\label{prg1}
Assume that the hypothesis $(H_{1}),(H_{2}) ~~\mbox{and}~~(H_{3})$ hold. Then,
 \begin{enumerate}
 \item[(i)] there exist  $\rho,~\beta_{0}>0$  such that $\mathcal{J}(u)\geq \beta_{0}$ for all
 $u\in E$ with $\|u\|=\rho$.
\item[(ii)] Let $\phi_{1} \in E\backslash \{0\}$. Then, $\mathcal{J}(t\phi_{1})\rightarrow -\infty,~~\mbox{as}~~t\rightarrow+\infty$.
 \end{enumerate}
\end{proposition}
 \textit{Proof}~~From (\ref{imp}), for all $\varepsilon>0$, there exists $C>0$ such that
\begin{equation*}
F(t)\leq \varepsilon |t|^{\frac{N}{2}}+C |t|^{q}\big(e^{a~t^{\gamma}}-1\big),~~~~~~\mbox{for all}~~t\in \R.
\end{equation*}

Then, using the last inequality, we get
$$\mathcal{J}(u)\geq \frac{2}{N}\|u\|^{\frac{N}{2}}-\varepsilon\int_{\mathbb{R}^{N}}|u|^{\frac{N}{2}}dx-C\int_{\mathbb{R}^{N}} |u|^{q}\big(e^{a~u^{\gamma}}-1\big)~dx.$$
 From the  H\"{o}lder inequality and using the following inequality
 $$\big(e^{s}-1\big)^{\nu}\leq e^{\nu s}-1,~~\forall~~s\geq0~~\forall\nu\geq1,$$ we obtain
\begin{equation}\label {eq:4.5}
\mathcal{J}(u)\geq \frac{2}{N}\|u\|^{\frac{N}{2}}-\varepsilon\int_{\mathbb{R}^{N}}|u|^{\frac{N}{2}}dx-C  (\int_{\mathbb{R}^{N}}\big(e^{2a~|u|^{\gamma}}-1\big)dx\Big)^{\frac{1}{2}}\|u\|^{q}_{2q}.
\end{equation}
From the Theorem \ref{th1.1}, if we choose $u\in E$ such that
\begin{equation}\label {eq:4.6}
2a \|u\|^{\gamma}\leq \alpha_{\beta},
\end{equation}
we get
$$\int_{\mathbb{R}^{N}}\big(e^{2a~|u|^{\gamma}}-1\big)dx=\int_{\mathbb{R}^{N}}\big(e^{2a~\|u\|^{\gamma}(\frac{|u|}{\|u\|})^{\gamma}})-1\big)dx<+\infty.$$
On the other hand from Sobolev embedding Lemma \ref{lem1}, there exist constants $C_{1}>0$ and $C_{2}>0$ such that  $\|u\|_{2q}\leq C_{2} \|u\|$ and   $\|u\|^{\frac{N}{2}}_{\frac{N}{2}}\leq C_{1} \|u\|^{\frac{N}{2}}$ .  So,
$$\mathcal{J}(u)\geq\frac{2}{N}\|u\|^{\frac{N}{2}}-\varepsilon C_{1}\|u\|^{\frac{N}{2}}-C\|u\|^{q}=\|u\|^{\frac{N}{2}}\big(\frac{2}{N}-\varepsilon C_{1}-C_{2}\|u\|^{q-\frac{N}{2}}\big),$$
for all $u\in E$ satisfying (\ref{eq:4.6}).
Since $q>N$, we can choose $\rho=\|u\|\leq (\frac{\alpha_{\beta}}{2a})^{\frac{1}{\gamma}}$ and for $\varepsilon$ such that  $ \displaystyle\frac{2}{NC_{1}}>\varepsilon $,  there exists $\beta_0=\rho^{\frac{N}{2}}\big(\displaystyle\frac{2}{N}-\varepsilon C_{1}-C_{2}\rho^{q-\frac{N}{2}}\big)>0$ with $\mathcal{J}(u) \geq\beta_0>0$. \hfill\\

(ii)  Let $\phi_{1}\in E\backslash\{0\}, $ $\|\phi_{1}\|=1$. We define the function
$$\displaystyle\varphi
(t)=\mathcal{J}(t\phi_{1})=\frac{2}{N}t^{\frac{N}{2}}\|\phi_{1}\|^{\frac{N}{2}}-\int_{\mathbb{R}^{4}}F(t\phi_{1})
dx.$$ By $(H_{2})$ and  $(H_{3})$ there exist two positive constants $C_{1}$  and $C_{2}$ such that
 $$F(t)\geq C_{1} t^{\theta}-C_{2}t^{\frac{N}{2}},~\forall~t\in \mathbb{R}.$$
 Hence, $$\varphi
(t)=\mathcal{J}(t\phi_{1})\leq\frac{2}{N}t^{\frac{N}{2}}\|\phi_{1}\|^{\frac{N}{2}}-C_{1}|t|^{\theta}\|\phi_{1}\|_{\theta}+C_{2}|t|^{\frac{N}{2}}\|\phi_{1}\|^{\frac{N}{2}}_{\frac{N}{2}}\rightarrow
-\infty,~~\mbox{as}~~t\rightarrow+\infty,$$ and it's easy to
conclude.\\\\

 \section{ Proof of Theorem \ref{th1.5} and Theorem \ref{th1.6}   }
\subsection{Palais-Smale sequence}
 We begin by the Palais-Smale sequence that is

 \begin{lemma}\label {cfF}Assume that $(H_{1}),(H_{2}),(H_{3})~~\mbox{and}~~(H_{4})$. If $(u_{n})\subset E$ is a (PS) sequence and $u\in E$ is a weak limit, then
 $$\int_{\mathbb{R}^{N}}F(u_{n})dx\rightarrow \int_{\mathbb{R}^{N}}F(u)dx.$$
 \end{lemma}
 Proof: We claim that \begin{align}\label{1}f(u_{n})\rightarrow f(u)~~\mbox{in}~~L^{1}(B_{R}),~~\mbox{for any}~~R>0.\end{align}
Since $u_{n}\rightarrow u~~\mbox{in}~~L^{\frac{N}{2}}(\mathbb{R}^{N}$), then $u_{n}\rightarrow u~~\mbox{in}~~L^{\frac{N}{2}}( B_{R})$. Furthermore,
$\int_{B_{R}}f(u_{n})u_{n}dx\leq  C$. It follows from Lemma \ref{lem2}  that (\ref{1}) holds.\\ Now from (\ref{1}) , we deduce that , for any $R>0$
\begin{align}\label{2} \int_{B_{R}}F(u_{n})dx\rightarrow \int_{B_{R}}F(u)dx.\end{align}
Indeed by $(H_{4})$ we have $$0<F(u_{n} )\leq M_{0}|f(u_{n})|~~\mbox{a.e}~~\{x\in \mathbb{R}^{N}~|~|u_{n}|\geq t_{0}\}$$ and from $(H_{2})$
$$0<F(u_{n} )\leq \frac{t_{0}}{\theta}|f(u_{n})|~~\mbox{a.e}~~\{x\in \mathbb{R}^{N}~|~|u_{n}|< t_{0}\}.$$
 Hence, applying the generalized Lebesgue dominated convergence theorem, we can conclude that (\ref{2}) holds for any $R > 0$.\\
 Now we claim that for any $\varepsilon>0$, there exists $R>1$ such that  \begin{align}\label{3} \bigg|\int_{\mathbb{R}^{N}\backslash B_{R}}(F(u_{n})-F(u))dx\bigg|\leq\varepsilon, \forall~n\geq 1.\end{align}
 In order to prove our claim, it's sufficient to see that for any $0<\varepsilon<1$, there exists $R>1$ such that
 \begin{align}\label{4}\int_{\mathbb{R}^{N}\backslash B_{R}}F(u_{n})dx\leq C\varepsilon~~\mbox{and}~~\int_{\mathbb{R}^{N}\backslash B_{R}}F(u)dx\leq C\varepsilon~\end{align}
 Let $R > 1$ arbitrarily fixed. By (\ref{eq:1.10}), and for $q=\frac{N}{2}$  there exists a positive constant $C_{3}$ such that
 \begin{align}\nonumber\int_{\mathbb{R}^{N}\backslash B_{R}}(F(u_{n})dx\leq \varepsilon\int_{\mathbb{R}^{N}\backslash B_{R}}|u_{n}|^{\frac{N}{2}}dx+C_{3}\int_{\mathbb{R}^{N}\backslash B_{R}}|u|^{\frac{N}{2}}_{n}~(e^{a~|u_{n}|^{\gamma}}-1)dx,~~\forall~n\geq 1, a>\alpha_{0}.\end{align}
 Using the power series expansion of the exponential function and estimating the
single terms with the radial lemma (12), the fact that $(u_{n})$ is bounded in $E$,  we get for any $n\geq 1$
\begin{align}\nonumber \int_{\mathbb{R}^{N}\backslash B_{R}}|u_{n}|^{\frac{N}{2}}~(e^{a~|u_{n}|^{\gamma}}-1)dx&=\sum_{k=1}^{+\infty}\frac{a^k}{k!} \int_{\mathbb{R}^{N}\backslash B_{R} } \vert u_{n}\vert^{\gamma k+\frac{N}{2}}dx\\&\nonumber\leq C_{4}\sum_{k=1}^{+\infty}\frac{a^k}{k!}(C)^{\gamma k}\|u_{n}\|^{\gamma k+\frac{N}{2}}\frac{R^{-2\frac{(N-1)(\gamma k)}{N}+1}}{\frac{2(N-1)(\gamma k)}{N}-1}\\&\leq C_{4}\frac{\|u_{n}\|^{\frac{N}{2}}}{R}\sum_{k=1}^{+\infty}\frac{a^k}{k!}(C\|u_{n}\|)^{\gamma k}\frac{1}{\frac{2(N-1)(\gamma )}{N}-1}\\&\nonumber\leq \frac{C_{4}}{R}\frac{\|u_{n}\|^{\frac{N}{2}}}{\frac{2(N-1)(\gamma )}{N}-1}\sum_{k=1}^{+\infty}\frac{(C~a~\|u_{n}\|^{\gamma)^{k}}}{k!}=\frac{C_{4}}{R}\frac{\|u_{n}\|^{\frac{N}{2}}}{\frac{2(N-1)(\gamma )}{N}-1}e^{a~C~\|u_{n}\|^{\gamma}},~~\forall n\\&\leq\frac{C_{5}(a,\varepsilon)}{R}\cdot\end{align}
On the other hand, using Sobolev embedding and the fact that  $(u_{n})$ is bounded in $E$, there exists $C_{1}$ such that
\begin{align} \nonumber \varepsilon\int_{\mathbb{R}^{N}\backslash B_{R}}|u_{n}|^{\frac{N}{2}}dx\leq\varepsilon C\|u_{n}\|^{\frac{N}{2}}\leq \varepsilon C_{1}.
\end{align}
It follows that, \begin{align}\nonumber \int_{\mathbb{R}^{N}\backslash B_{R}}(F(u_{n})dx\leq\varepsilon C_{1}+\frac{C_{5}(a,\varepsilon)}{R}\cdot\end{align}
We can assume without loss of generality that $\frac{C_{5}(a,\varepsilon)}{\varepsilon}>1$. Taking $$R=\frac{C_{5}(a,\varepsilon)}{\varepsilon}>1,$$ we get
\begin{align}\nonumber\int_{\mathbb{R}^{N}\backslash B_{R}}F(u_{n})dx\leq\varepsilon (C_{1}+1),~~\forall~~n\geq1.\end{align}
In the same way,\begin{align}\nonumber\int_{\mathbb{R}^{N}\backslash B_{R}}F(u)dx\leq\varepsilon (C_{1}+1).\end{align}
Then,
\begin{align}\nonumber\bigg|\int_{\mathbb{R}^{N}\backslash B_{R}}(F(u_{n})-F(u))dx\bigg|\leq \int_{\mathbb{R}^{N}\backslash B_{R}}F(u_{n})+\int_{\mathbb{R}^{N}\backslash B_{R}}F(u)dx\leq 2\varepsilon (C_{1}+1)\end{align}
and (\ref{3}) holds.
\subsection{Estimate of the mountain pass level}


\begin{lemma}\label{lem4.2}Assume that $f$ verifies the conditions $(H_{1}),(H_{2}),(H_{4})~~\mbox{and}~(H_{5})$. Then, for the sequence $(v_{n})$ given by (\ref{eq:5.2}), there exists $n\geq1$ such that
  \begin{equation}\label{eq:5.3}
 \max_{t\geq 0}\mathcal{J}(tv_{n})<\displaystyle\frac{2}{N}(\frac{\alpha_{\beta}}{\alpha_{0}})^{\frac{N}{2\gamma}}\cdot
       \end{equation}
\end{lemma}
{\it Proof}~~ By contradiction, suppose that for all $n\geq1$,
$$ \max_{t\geq 0}\mathcal{J}(tv_{n})\geq\displaystyle\frac{2}{N}(\frac{\alpha_{\beta}}{\alpha_{0}})^{\frac{N}{2\gamma}}\cdot$$
 By contradiction, suppose that for all $n\geq1$,
$$ \max_{t\geq 0}\mathcal{J}(tv_{n})\geq\displaystyle\frac{2}{N}(\frac{\alpha_{\beta}}{\alpha_{0}})^{\frac{N}{2\gamma}}\cdot$$
Therefore,
for any $n\geq1$, there exists $t_{n}>0$ such that
$$\max_{t\geq0}\mathcal{J}(tv_{n})=\mathcal{J}(t_{n}v_{n})\geq \displaystyle\frac{2}{N}(\frac{\alpha_{\beta}}{\alpha_{0}})^{\frac{N}{2\gamma}}$$
and so,
$$\frac{2}{N}t_{n}^{\frac{N}{2}}-\int_{B}F(x,t_{n}v_{n})dx\geq \displaystyle\frac{2}{N}(\frac{\alpha_{\beta}}{\alpha_{0}})^{\frac{N}{2\gamma}}\cdot$$
Then, by using $(H_{1})$
   \begin{equation}\label{eq:5.4}
t_{n}^{\frac{N}{2}}\geq\displaystyle(\frac{\alpha_{\beta}}{\alpha_{0}})^{\frac{N}{2\gamma}}\cdot
       \end{equation}
On the other hand,
 $$\frac{d}{dt}\mathcal{J}(tv_{n})\big|_{t=t_{n}}=t^{\frac{N}{2}-1}_{n}-\int_{B}f(x,t_{n}v_{n})v_{n}dx=0,$$
then,
\begin{equation}\label{eq:5.5}
t_{n}^{\frac{N}{2}}=\int_{B}f(x,t_{n}v_{n})t_{n}v_{n}dx.
\end{equation}

Now, we claim that the sequence $(t_{n})$ is bounded in $(0,+\infty)$. Indeed,
 it follows from $(H_{4})$ that for all $\varepsilon>0$, there exists
$t_{\varepsilon}>0$ such that
\begin{equation}\label{eq:5.6}
f(t)t\geq
(\gamma_{0}-\varepsilon)e^{\alpha_{0}t^{\gamma}}~~\forall
|t|\geq t_{\varepsilon}.
\end{equation}
Using  (\ref{eq:5.5}), we get
\begin{equation*}\label{eq:5.13}
t_{n}^{\frac{N}{2}}=\int_{B}f(x,t_{n}v_{n})t_{n}v_{n}dx\geq
\int_{0\leq |x|\leq \frac{1}{\sqrt[N]{n}}}f(x,t_{n}v_{n})t_{n}v_{n}dx
\cdot\end{equation*}
We have  for all $0\leq |x|\leq \frac{1}{\sqrt[4]{n}}$, $\displaystyle w^{\gamma}_{n}\geq \displaystyle\bigg(\frac{\log (e\sqrt[N]{n} )}{\alpha_{\beta}}\bigg)\cdot$
From (\ref{eq:5.4}) and the result of Lemma \ref{lem6},  $$t_{n} v_{n}\geq\frac{t_{n}}{\|w_{n}\|}\big(\frac{\log (e\sqrt[N]{n} )}{\alpha_{\beta}}\big)^{\frac{1}{\gamma}} \rightarrow\infty~~\mbox{as}~~n\rightarrow+\infty.$$Hence,
it follows from  (\ref{eq:5.6}) that for all $\varepsilon >0$, there exists $n_{0}$ such that for all $n\geq n_{0}$
\begin{equation*}
t^{\frac{N}{2}}_{n}\geq \displaystyle (\gamma_{0}-\varepsilon)\int_{0\leq |x|\leq \frac{1}{\sqrt[N]{n}}}e^{\alpha_{0}t^{\gamma}_{n}|v|^{\gamma}_{n}}dx
  \end{equation*}
and  \begin{equation}\label{eq:5.7}
 t^{\frac{N}{2}}_{n}\geq \displaystyle NV_{N} (\gamma_{0}-\varepsilon)\int_{0}^{\frac{1}{\sqrt[N]{n}}}r^{N-1}e^{\alpha_{0}t^{\gamma}_{n} \displaystyle\big(\frac{\log (e\sqrt[N]{n})}{\|w_{n}\|^{\gamma}\alpha_{\beta}}\big)}dr\cdot
\end{equation}
Hence,
 \begin{equation*}\label{eq:5.13}
1 \geq  \displaystyle NV_{N}  (\gamma_{0}-\varepsilon)~~\displaystyle e^{\alpha_{0}t^{\gamma}_{n} \displaystyle\big(\frac{\log (e\sqrt[N]{n})}{\|w_{n}\|^{\gamma}\alpha_{\beta}}\big)-\log Nn-\frac{N}{2}\log t_{n}}.
 \end{equation*}
 Therefore $(t_{n})$ is bounded. Also, we have from the formula  (\ref{eq:5.5}) that
$$\displaystyle\lim_{n\rightarrow+\infty} t_{n}^{\frac{N}{2}}\geq\displaystyle(\frac{\alpha_{\beta}}{\alpha_{0}})^{\frac{N}{2\gamma}}\cdot$$
Now, suppose that
$$\displaystyle\lim_{n\rightarrow+\infty} t_{n}^{\frac{N}{2}}>\displaystyle(\frac{\alpha_{\beta}}{\alpha_{0}})^{\frac{N}{2\gamma}},$$
then for $n$ large enough, there exists some $\delta>0$ such that $ t_{n}^{\gamma}\geq \frac{\alpha_{\beta}}{\alpha_{0}}+\delta$. Consequently the right hand side of (\ref{eq:5.7}) tends to infinity and this contradicts the boudness of  $(t_{n})$. Since $(t_{n})$ is bounded, we get
\begin{equation}\label{eq:5.8}
 \displaystyle
 \lim_{n\rightarrow+\infty}t_{n}^{\frac{N}{2}}=\displaystyle (\frac{\alpha_{\beta}}{\alpha_{0}})^{\frac{N}{2\gamma}}\cdot
 \end{equation}

Let consider the unit ball $B$ of $\mathbb{R}^{N}$ and  the sets $$\mathcal{A}_{n}=\{x\in B~|~ t_{n}v_{n}\geq
t_{\varepsilon}\}~~\mbox{and}~~\mathcal{C}_{n}=B\setminus \mathcal{A}_{n}.$$ we have,
$$\begin{array}{rclll} t_{n}^{\frac{N}{2}}&=&\displaystyle\int_{\mathbb{R}^{N}}f(t_{n}v_{n})t_{n}v_{n}dx\geq\displaystyle\int_{B}f(t_{n}v_{n})t_{n}v_{n}dx=\int_{\mathcal{A}_{n}}f(t_{n}v_{n})t_{n}v_{n}dx+\int_{\mathcal{C}_{n}}f(t_{n}v_{n})t_{n}v_{n} $$\\
&\geq& \displaystyle(\gamma_{0}-\varepsilon)\int_{\mathcal{A}_{n}}e^{\alpha_{0}t_{n}^{\gamma}v_{n}^{\gamma}}dx + \int_{\mathcal{C}_{n}}f(t_{n}v_{n})t_{n}v_{n}dx\\
&=&\displaystyle(\gamma_{0}-\varepsilon)\int_{B}e^{\alpha_{0}t_{n}^{\gamma}v_{n}^{\gamma}}dx-
(\gamma_{0}-\varepsilon)\int_{\mathcal{C}_{n}}e^{\alpha_{0}t_{n}^{\gamma}v_{n}^{\gamma}}dx\\ &+&\displaystyle\int_{\mathcal{C}_{n}}f(t_{n}v_{n})t_{n}v_{n}dx.
\end{array}$$
Since $v_{n}\rightarrow 0 ~~\mbox{a.e in }~~B$,
$\chi_{\mathcal{C}_{n}}\rightarrow1~~\mbox{a.e in}~~B$, therefore using the dominated convergence
theorem, we get $$\displaystyle\int_{\mathcal{C}_{n}}f(x,t_{n}v_{n})t_{n}v_{n}dx\rightarrow 0~~\mbox{and}~~\int_{\mathcal{C}_{n}}e^{\alpha_{0}t_{n}^{\gamma}v_{n}^{\gamma}}dx\rightarrow NV_{N}\cdot$$Then,$$\ \lim_{n\rightarrow+\infty} t_{n}^{\frac{N}{2}}=\displaystyle(\frac{\alpha_{\beta}}{\alpha_{0}})^{\frac{N}{2\gamma}}\geq(\gamma_{0}-
\varepsilon)\lim_{n\rightarrow+\infty}\int_{B}e^{\alpha_{0}t_{n}^{\gamma}v_{n}^{\gamma}}dx-(\gamma_{0}-
\varepsilon)NV_{N}\cdot$$
On the other hand,
$$\int_{B}e^{\alpha_{0}t_{n}^{\gamma}v_{n}^{\gamma}}dx \geq \int_{\frac{1}{\sqrt[N]{n}}\leq |x|\leq \frac{1}{2}}e^{\alpha_{0}t_{n}^{\gamma}v_{n}^{\gamma}}dx+\int_{\mathcal{C}_{n}} e^{\alpha_{0}t_{n}^{\gamma}v_{n}^{\gamma}}dx\cdot$$

Then, using (\ref{eq:5.4})  $$\lim_{n\rightarrow+\infty} t_{n}^{\frac{N}{2}}\geq \lim_{n\rightarrow+\infty}\displaystyle(\gamma_{0}-\varepsilon)\int_{B}e^{\alpha_{0}t_{n}^{\gamma}v_{n}^{\gamma}}dx\geq \displaystyle\lim_{n\rightarrow+\infty}(\gamma_{0}-\varepsilon) NV_{N}\int^{\frac{1}{2}}_{\frac{1}{\sqrt[N]{n}}}r^{N-1} e^{C^{\gamma}(N,\beta)\frac{ \big(\log \frac{e}{r}\big)^{\frac{N}{N-2}}}{(\log (e\sqrt[N]{n}))^{\frac{2}{N-2}}\|w_{n}\|^{\gamma}}} dr.$$

Therefore,  making the change of variable$$s=\frac{C(N,\beta)^{\gamma}(\log \frac{e}{r})}{(\log (e\sqrt[N]{n}))^{\frac{2}{N-2}}\|w_{n}\|^{\gamma}}=P\frac{(\log \frac{e}{r})}{\|w_{n}\|^{\gamma}},~~\mbox{with}~~P=\frac{C(N,\beta)^{\gamma}}{(\log (e\sqrt[N]{n}))^{\frac{2}{N-2}}}$$
we get
$$\begin{array}{rclll}\displaystyle \lim_{n\rightarrow+\infty} t_{n}^{\frac{N}{2}}&\geq & \displaystyle\lim_{n\rightarrow+\infty}\displaystyle(\gamma_{0}-\varepsilon)\int_{B}e^{\alpha_{0}t_{n}^{\gamma}v_{n}^{\gamma}}dx\\\\
&\geq&\displaystyle \lim_{n\rightarrow+\infty} NV_{N}(\gamma_{0}-\varepsilon)\frac{\|w_{n}\|^{\gamma}}{P}\int^{\frac{P \log(e\sqrt[N]{n})}{\|w_{n}\|^{\gamma}}}_{\frac{P \log(2e)}{\|w_{n}\|^{\gamma}}}\displaystyle e^{N(1-\frac{s\|w_{n}\|^{\gamma}}{P})+\frac{\|w_{n}\|^{\frac{2\gamma}{N-2}}}{P^{\frac{N}{N-2}}}s^{\frac{N}{N-2}}}ds\\\\
&\geq&\displaystyle \lim_{n\rightarrow+\infty} NV_{N}(\gamma_{0}-\varepsilon)\frac{\|w_{n}\|^{\gamma}}{P}e^{N}\int^{\frac{P \log(e\sqrt[N]{n})}{\|w_{n}\|^{\gamma}}}_{\frac{P \log(2e)}{\|w_{n}\|^{\gamma}}}\displaystyle e^{-\frac{N}{P}\|w_{n}\|^{\gamma}s}ds\\
&=&\displaystyle\lim_{n\rightarrow+\infty}(\gamma_{0}-\varepsilon)NV_{N}\frac{e^{N}}{N}\big(-e^{-N\log (e\sqrt[N]{n})}+e^{-N\log (2e)}\big)\\
&=&\displaystyle(\gamma_{0}-\varepsilon)V_{N}e^{N(1-\log (2e))}\cdot
\end{array}$$

It follows that
\begin{equation*}\label{eq:4.15}
\displaystyle\displaystyle(\frac{\alpha_{\beta}}{\alpha_{0}})^{\frac{N}{2\gamma}}\geq (\gamma_{0}-\varepsilon)V_{N}e^{N(1-\log (2e))}
\end{equation*}

for all $\varepsilon>0$. So,
$$\gamma_{0}\leq\displaystyle \frac{(\frac{\alpha_{\beta}}{\alpha_{0}})^{\frac{N}{2\gamma}}}{V_{N}e^{N(1-\log (2e))}}, $$
   which is in contradiction with  the condition $(H_{5})$.
   \subsection{The compactness level of the energy}

The primary challenge within the variational approach to the critical growth problem arises due to the absence of compactness. Specifically, the global Palais-Smale condition doesn't hold. However, a partial Palais-Smale condition is retained under a specific threshold. In the subsequent proposition, we pinpoint the initial level at which the energy exhibits non-compactness.
   \begin{proposition}\label{propPS} Let $\mathcal{J}$ be the energy associated to the problem (\ref{eq:1.1}) defined by (\ref{energy}), and suppose that the conditions  $(H_{1})$, $(H_{2})$ and  $(H_{4})$ are satisfied.
    \begin{itemize}\item [(i)]If the function $f(t)$ satisfies the condition (\ref{eq:1.4}) for some
$\alpha_{0} >0$, then the functional $\mathcal{J}$ satisfies the Palais-Smale condition $(PS)_{c}$ for any
$$c<\displaystyle\frac{2}{N}(\frac{\alpha_{\beta}}{\alpha_{0}})^{\frac{N}{2\gamma}}.$$

\item  [(ii)]If $f$ is subcritical at $+\infty$, then the functional $\mathcal{J}$ satisfies the Palais-Smale condition $(PS)_{c}$ for any $c\in \mathbb{R}$.
    \end{itemize}

\end{proposition}
 {\it Proof}~~$(i)$~~Consider a $(PS)_{c}$ sequence ($u_{n}$) in $E$, for some $c\in \R$, that is
\begin{equation}\label{eq:4.7}
\mathcal{J}(u_{n})=\frac{2}{N}\|u_{n}\|^{\frac{N}{2}}-\int_{B}F(x,u_{n})dx \rightarrow c ,~~n\rightarrow +\infty
\end{equation}
and
 \begin{align}\label{eq:4.8}
\mathcal{ J}'(u_{n})\varphi =& \nonumber \Big|\int_{B}\big( w_{\beta}(x)~|\Delta u_{n}|^{\frac{N}{2}-2}~\Delta u_{n}~~\Delta\varphi+|\nabla u_{n}|^{\frac{N}{2}-2} \nabla u_{n}.\nabla \varphi +V(x)|u_{n}|^{\frac{N}{2}-2}u_{n}\varphi\big)~dx-\int_{B}f(x,u_{n})~ \varphi~dx,
\Big|\\&\leq \varepsilon_{n}\|\varphi\|,
\end{align}
for all $\varphi \in E$, where $\varepsilon_{n}\rightarrow0$, as $n\rightarrow +\infty$.\\
For $n$ large enough, there exists a constant $C>0$ such that
$$\displaystyle\frac{2}{N}\|u_{n}\|^{\frac{N}{2}}\leq C+\displaystyle\int_{B}F(x,u_{n})dx,$$
 From $(H_{2})$, it follows that $$ \int_{\mathbb{R}^{N}}F(u_{n})dx\leq \frac{1}{\theta} \int_{\mathbb{R}^{N}}f(u_{n})u_{n}dx.$$
 Using (\ref{eq:4.8}) with $\varphi=u_{n}$, we obtain
 $$\int_{\mathbb{R}^{N}}f(u_{n})u_{n}dx\leq \varepsilon_{n}\|u_{n}\|+\|u_{n}\|^{\frac{N}{2}}.$$
 Therefore,

$$\displaystyle\frac{2}{N}\|u_{n}\|^{\frac{N}{2}}\leq C_{1}+ \frac{\varepsilon_{n}}{\theta}\|u_{n}\|+ \frac{1}{\theta} \|u_{n}\|^{\frac{N}{2}},$$

Since, $\theta>\frac{N}{2}$, we get

$$0<\displaystyle(\frac{2}{N}-\frac{1}{\theta})\|u_{n}\|^{\frac{N}{2}}\leq C+ \frac{\varepsilon_{n}}{\theta}\|u_{n}\|.$$

 We deduce that the sequence $(u_{n})$ is bounded in $E$. As consequence, there exists $u\in E$ such that, up to subsequence,
 $u_{n}\rightharpoonup u $ weakly in $E$, $u_{n}\rightarrow u$ strongly in $L^{q}(B)$, for all $1\leq q\leq N$ and $u_{n}(x)\rightarrow u(x)$ a.e. in $\mathbb{R}^{N}$.  Also, we can follow  \cite{CJ} to prove that  $\nabla u_{n}(x)\rightarrow\nabla u(x) ~\mbox{a.e}~x\in\mathbb{R}^{N} $ and  $\Delta u_{n}(x)\rightarrow\Delta u(x) ~\mbox{a.e}~x\in\mathbb{R}^{N} .$\\
Furthermore, we have, from (\ref{eq:4.7}) and (\ref{eq:4.8}), that
\begin{equation}\label{eq:4.11}
0<\int_{\mathbb{R}^{N}} f(x,u_{n})u_{n}\leq C,
 \end{equation}
and
 \begin{equation}\label{eq:4.12}
0<\int_{\mathbb{R}^{N}} F(x,u_{n})\leq C.
 \end{equation}
 By Lemma \ref{cfF} , we have

\begin{equation}\label{eq:4.14}
F(x,u_{n})\rightarrow F(x,u) ~~\mbox{in}~~L^{1}(\mathbb{R}^{N}) ~~as~~ n\rightarrow +\infty.
 \end{equation}
 Then, from (\ref{eq:4.7}), we get

\begin{equation}\label{eq:4.15}
\displaystyle\lim_{n\rightarrow+\infty}\int_{B}f(x,u_{n})u_{n}dx=\frac{N}{2}(c+\int_{B}F(x,u)dx)
 \end{equation}
and from (\ref{eq:4.8}), we have
\begin{equation}\label{eq:4.16}
\displaystyle\lim_{n\rightarrow+\infty}\int_{B}f(x,u_{n})u_{n}dx=\frac{N}{2}(c+\int_{B}F(x,u)dx).
 \end{equation}
It follows from $(H_{2})$ and (\ref{eq:4.8}), that
\begin{equation}\label{eq:4.17}
\displaystyle
\lim_{n\rightarrow+\infty}\frac{N}{2} \int_{B}F(x,u_{n})dx\leq \displaystyle
\lim_{n\rightarrow+\infty}\int_{B}f(x,u_{n})u_{n}dx= \frac{N}{2}(c+\int_{B}F(x,u)dx).
\end{equation}
Then, passing to the limit in (\ref{eq:4.8}) and using (\ref{eq:4.16}),we obtain  that $u$ is a weak solution of the problem (\ref{eq:1.1}) that is $$\int_{B}\big( w_{\beta}(x)~|\Delta u|^{\frac{N}{2}-2}~\Delta u~~\Delta\varphi+|\nabla u|^{\frac{N}{2}-2} \nabla u.\nabla \varphi +V(x)|u|^{\frac{N}{2}-2}u\varphi\big)~dx=\int_{B}f(x,u)~ \varphi~dx,~~~~~~\mbox{for all }~~\varphi \in E.$$ Taking $\varphi=u$ as a test function, we get
$$\|u\|^{\frac{N}{2}}=\int_{B} w_{\beta}(x) |\Delta u|^{\frac{N}{2}} dx+\int_{B}|\nabla u|^{\frac{N}{2}}dx+\int_{B}V(x)|u|^{\frac{N}{2}}dx=
\int_{B}f(x,u)u dx\geq \frac{N}{2}\int_{B}F(x,u)dx\cdot$$
~Hence $\mathcal{J}(u)\geq 0$ . We also have by the Fatou's lemma and (\ref{eq:4.14})  $$0\leq \mathcal{J}(u)\leq \frac{2}{N}\liminf_{n\rightarrow\infty} \|u_{n}\|^{\frac{N}{2}}-\int_{B}F(x,u)dx=c.$$

So, we will finish the proof by considering  three cases for the level $c$.\\\\
{\it Case 1.} $c=0$.~~ In this case
 $$0\leq\mathcal{ J}(u)\leq \liminf_{n\rightarrow+\infty}\mathcal{J}(u_{n})=0.$$
So, $$\mathcal{J}(u)=0$$
   and then by (\ref{eq:4.14})
$$\displaystyle\lim_{n\rightarrow +\infty}\frac{2}{N}\|u_{n}\|^{\frac{N}{2}}=\int_{B}F(x,u)dx=\frac{2}{N}\|u\|^{\frac{N}{2}}.$$
By Brezis-Lieb's Lemma \cite{Br}, it follows that $u_{n}\rightarrow u~~\mbox{in}~~E$.\\
\noindent {\it Case 2.} $c>0$ and $u=0$. We will prove that this is not possible.\\
We will show in the following that \begin{align}\label{eq:4.18}\mathcal{J}_{n}=\int_{\mathbb{R}^{N}}f(u_{n})u_{n}dx\rightarrow0~~\mbox{as}~~n\rightarrow+\infty\cdot\end{align}
In fact, if (\ref{eq:4.18}) holds, then taking $\varphi=u_{n}$ in (\ref{eq:4.8}) we get
$$\|u_{n}\|^{\frac{N}{2}}\leq c\varepsilon_{n}+\mathcal{J}_{n}\rightarrow0~~\mbox{as}~~n\rightarrow0$$
which gives $\displaystyle\lim_{n\rightarrow+\infty}\|u_{n}\|^{\frac{N}{2}}=0$. This contradicts the fact that $c\neq0$ because from (\ref{eq:4.7}) and (\ref{eq:4.8}), we have

$$\displaystyle\lim_{n\rightarrow +\infty}\|u_{n}\|^{\frac{N}{2}}=\frac{N}{2}c.$$
It therefore remains to prove that (\ref{eq:4.18}) is valid. Let $a>\alpha_{0}$ and $q\geq 1$. By (\ref{imp}) and again the boundedness of $(u_{n})$ in $E$, we obtain for all $\varepsilon>0$,
\begin{align}\nonumber \mathcal{J}_{n}\leq C_{1}\varepsilon +C(a,q,\varepsilon)\int_{\mathbb{R}^{N} }|u|^{q}_{n}~(e^{a~|u_{n}|^{\gamma}}-1)dx~~\forall~~n\geq1.\end{align}
Applying H\"{o}lder's inequality with $p,p'>1$ and $\frac{1}{p}+\frac{1}{p'}=1$, we get
\begin{align}\nonumber \mathcal{ J}_{n}(u)\leq  C\|u_{n}\|^{q}_{p'q}\bigg(\int_{\mathbb{R}^{N} }~(e^{p~a~|u_{n}|^{\gamma}}-1)dx\bigg)^{\frac{1}{p}}\end{align}

Since $(\frac{N}{2}c)^{\frac{2\gamma}{N}}<\displaystyle(\frac{\alpha_{\beta}}{\alpha_{0}})$, there exists
$\eta\in(0,\frac{1}{2})$ such that
$(\frac{N}{2}c)^{\frac{2\gamma}{N}}=\displaystyle(1-2\eta)\displaystyle(\frac{\alpha_{\beta}}{\alpha_{0}})$.
On the other hand, $\|u_{n}\|^{\gamma}\rightarrow
(\frac{N}{2}c)^{\frac{2\gamma}{N}}$, so there exists $n_{\eta}>0$ such that for
all $n\geq n_{\eta}$, we get $\|u_{n}\|^{\gamma}\leq
(1-\eta)\frac{\alpha_{\beta}}{\alpha_{0}}$.
Therefore, if we choose $a=(1+\frac{\eta}{2})\alpha_{0}$, $p=(1+\frac{\eta}{2})$ we get
$$p~a~(\frac{|u_{n}|}{\|u_{n}\|})^{\gamma}\|u_{n}\|^{\gamma}\leq\alpha_{0}(1+\frac{\eta}{2})^{2}(\frac{|u_{n}|}{\|u_{n}\|})^{\gamma}(1-\eta)\leq
~\alpha_{\beta}(\frac{|u_{n}|}{\|u_{n}\|})^{\gamma}\cdot$$

 Therefore, the  integral is
 bounded in view of (\ref{eq:1.4}). On the other hand, choosing $q > 2$ , so $p'q>2$ and therefore $u_{n}\rightarrow 0~~ L^{qp'}(\mathbb{R}^N)$. Then $\mathcal{J}_{n}\rightarrow0~~\mbox{as}~~n\rightarrow+\infty.$\\\\
{\it Case 3.} $c>0$ and $u\neq 0$.~~ In this case,
we claim that $\mathcal{J}(u)=c$ and therefore, we get
$$\lim_{n\rightarrow
+\infty}\|u_{n}\|^{\frac{N}{2}}=\frac{N}{2}\big(c+\int_{B}F(x,u)dx\big)=\big(\mathcal{J}(u)+\int_{B}F(x,u)dx\big)=\|u\|^{\frac{N}{2}}.$$
Do not forgot that
 $$\mathcal{J}(u)\leq \frac{2}{N}\liminf_{n\rightarrow+\infty} \|u_{n}\|^{\frac{N}{2}}-\int_{B}F(x,u)dx=c.$$
We argue by contradiction and suppose that $\mathcal{J}(u)<c$. Then,
\begin{equation}\label{eq:4.20}
\|u\|^{\frac{N}{2}}<(\frac{N}{2}\big(c+\int_{B}F(x,u)dx\big)\big)^{\frac{2}{N}}.
\end{equation}
Set
$$v_{n}=\displaystyle\frac{u_{n}}{\|u_{n}\|}$$ and
$$v=\displaystyle\frac{u}{(\frac{N}{2}\big(c+\displaystyle\int_{B}F(x,u)dx\big))^{\frac{2}{N}}}\cdot$$
We have $\|v_{n}\|=1$, $v_{n}\rightharpoonup v$  in $E$,  $\nabla v_{n}(x)\rightarrow\nabla v(x) ~\mbox{a.e}~x\in\mathbb{R}^{N} ,$  $\Delta v_{n}(x)\rightarrow\Delta v(x) ~\mbox{a.e}~x\in\mathbb{R}^{N} $ $v\not\equiv 0$ and $\|v\|<1$. So, by Lemma \ref{Lionstype}, we get
\begin{align}\label{sup}\displaystyle \sup_{n}\int_{\mathbb{R}^{N}}\big(e^{p \alpha_{\beta}|v_{n}|^{\gamma}}-1\big)dx<\infty,~\mbox{ for}~~1<p<U(v)=(1-\|v\|^{\frac{N}{2}})^{\frac{-2\gamma}{N}}.\end{align}
Since $u_{n} \hookrightarrow u~~\mbox{in}~~ E,$ it suffice to prove that
$$\mathcal{ J}'(u_{n})(u_{n}-u)\rightarrow 0~~\mbox{as}~n\rightarrow +\infty,$$
and that's the case when \begin{align}\label{eq:4.21}\int_{\mathbb{R}^{N}}f(u_{n})(u_{n}-u)dx\rightarrow 0.\end{align}

Arguing as in Case 1, we can thus reduce the proof of (\ref{eq:4.21}) to showing the existence of $a > \alpha_{0}$ and $q \geq1$ such that \begin{align}\nonumber \mathcal{I}_{n}:=\int_{\mathbb{R}^{N} }~|u_{n}|^{q-1}|u_{n}-u|(e^{a~|u_{n}|^{\gamma}}-1)dx\rightarrow0~~\mbox{as}~~n\rightarrow+\infty\cdot\end{align}
We apply H\"{o}lder's inequality twice with $p,p',t,t'>1$ and $\frac{1}{p}+\frac{1}{p'}=\frac{1}{t}+\frac{1}{t'}=1$, we get

\begin{align}\nonumber \mathcal{I}_{n}\leq C(a,\varepsilon)\|u_{n}\|^{q-1}_{p't'(q-1)}\|u_{n}-u\|_{p't}\bigg(\int_{\mathbb{R}^{N} }~(e^{p~a~|u_{n}|^{\gamma}}-1)dx\bigg)^{\frac{1}{p}}\leq\bigg(\int_{\mathbb{R}^{N} }~(e^{\tau~a~|u_{n}|^{\gamma}}-1)dx\bigg)^{\frac{1}{p}}\end{align} for any $\tau>p$.\newline
From (\ref{sup}), it follows that \begin{align}\nonumber \sup_{n}\int_{\mathbb{R}^{N}}\big(e^{\tau~a~ |u_{n}|^{\gamma}}-1\big)dx= \sup_{n}\int_{\mathbb{R}^{N}}\big(e^{\tau~a~ |v_{n}|^{\gamma}\|u_{n}\|^{\gamma}}-1\big)dx < \infty\end{align}
provided $a>\alpha_{0}, p>1 $ and $\tau>p$ can be chosen so that
$a~\tau\|u_{n}\|^{\gamma}< ~(1-\|v\|^{\frac{N}{2}})^{\frac{-2\gamma}{N}}
\alpha_{\beta}$ and $1<p<U(v)=(1-\|v\|^{\frac{N}{2}})^{\frac{-2\gamma}{N}}$.\\
We have
\begin{align}\label{p}\displaystyle(1-\|v\|^{\frac{N}{2}})^{\frac{-2\gamma}{N}}=\displaystyle\big(\frac{(\frac{N}{2}(c+\int_{B}F(x,u)dx)}{(\frac{N}{2}(c+\int_{B}F(x,u)dx))-\|u\|^{\frac{N}{2}})}\big)^{\frac{2\gamma}{N}}=
\big(\frac{c+\int_{B}F(x,u)dx}{c-\mathcal{J}(u)}\big)^{\frac{2\gamma}{N}}.\end{align}
On the other hand, \begin{align}\nonumber \displaystyle\lim_{n\rightarrow+\infty}\|u_{n}\|^{\gamma}=(\frac{N}{2}\big(c+\int_{B}F(x,u)dx)\big)^{\frac{2\gamma}{N}},\end{align}
then, for all $\varepsilon$ such that $0<\varepsilon<1$  and for $n$ large enough
\begin{align}\nonumber \alpha_{0}(1+\varepsilon)\|u_{n}\|^{\gamma}\leq \alpha_{0}(1+2\varepsilon)
(\frac{N}{2}\big(c+\int_{B}F(x,u)dx\big)^{\frac{2\gamma}{N}}.\end{align}

Taking $a=(1+\varepsilon)\alpha_{0}$ , $\tau=(1+\varepsilon)p$ and using (\ref{p}), we get
\begin{align}\nonumber a\|u_{n}\|^{\gamma}\tau\leq\alpha_{0}(1+7\varepsilon)(\frac{N}{2}\big(c+\int_{\mathbb{R}^{N}}F(x,u)dx\big)^{\frac{2\gamma}{N}}(1-\|v\|^{\frac{N}{2}})^{\frac{-2\gamma}{N}}\leq &  p\alpha_{0}(1+7\varepsilon)2^{\frac{2\gamma}{N}}\big(c-\mathcal{J}(u)\big)^{\frac{2\gamma}{N}}.\end{align}
But $\mathcal{J}(u)\geq 0$ and  $c<\displaystyle\frac{2}{N}(\frac{\alpha_{\beta}}{\alpha_{0}})^{\frac{N}{2\gamma}}$ , then there exists $\eta\in(0,1)$ such that $c^{\frac{2\gamma}{N}}=(1-\eta)(\frac{2}{N})^{\frac{2\gamma}{N}}\frac{\alpha_{\beta}}{\alpha_{0}}$.
\newline If we choose $\varepsilon=\frac{\eta}{7}$, we get,
\begin{align} \nonumber a\|u_{n}\|^{\gamma}\tau\leq (1+\eta)(1-\eta)p \alpha_{\beta}\leq p \alpha_{\beta}<p~~(1-\|v\|^{\frac{N}{2}})^{\frac{-2\gamma}{N}} .\end{align}
So, with this choice of $\tau>p>1$ and $a>\alpha_{0}$, we have
\begin{align}\nonumber \mathcal{I}_{n}\leq  C(a,\alpha_{0})\|u_{n}\|^{q-1}_{p't'(q-1)}\|u_{n}-u\|_{p't}
\end{align}
where $ C(a,\alpha_{0})$ is a positive constant depending only on $a ~\mbox{and}~\alpha_{0}.$ Now, since $(q-1)p't'>q-1$ and
$p't>t$, choosing $q\geq3$ and $t\geq2$ we have that $(u_{n})$ is bounded in $L^{((q-1)p't'}(\mathbb{R}^{N})$, so
$\mathcal{I}_{n}\rightarrow0~~\mbox{as}~~n\rightarrow+\infty$.


Hence,
$$\displaystyle\lim_{n\rightarrow+\infty}\|u_{n}\|^{\frac{N}{2}}=\frac{N}{2}(c+\int_{\mathbb{R}^{N}}F(x,u)dx)=\|u\|^{\frac{N}{2}}$$
and this contradicts (\ref{eq:4.20}). So, $\mathcal{J}(u)=c$ and  consequently, $u_{n} \rightarrow u$.\\
 {\it Proof}~~$(ii)$ Follows from $(i)$. Indeed, since $(u_{n} )$ is bounded in $E$, there exists a positive constant $M>0$ such that $\|u_{n}\|\leq M$. As $f$ is subcritical at infinity, then if we choose  $a>0$ such that
 $a\leq \frac{\alpha_{\beta}}{p~M^{\gamma}}$, the integral $$\int_{\mathbb{R}^{N} }~(e^{p~a~|u_{n}|^{\gamma}}-1)dx$$ is finite for all $p\geq1$. So, arguing as in $(i)$ we get that  \begin{align}\nonumber \mathcal{I}_{n}\rightarrow0~~\mbox{as}~~n\rightarrow+\infty.\end{align}
 \\Now according to the  Proposition \ref{propPS}, the functional $\mathcal{J}$  satisfies the $(PS)_{c}$ condition at a level $c<\displaystyle\frac{2}{N}(\frac{\alpha_{\beta}}{\alpha_{0}})^{\frac{N}{2\gamma}}$, in the critical case and at all level $c$, in the subcritical case . Moreover, Proposition \ref{prg1} confirms that the functional $\mathcal{J}$ exhibits a mountain pass structure. Consequently, by the Ambrosetti and Rabinowitz Theorem \cite{AR}, $\mathcal{J}$ possesses a non-zero critical point $u$ within the space $E$. This leads to the proof of Theorem \ref{th1.5} and Theorem \ref{th1.6}.

    \section*{ Statements and Declarations:}

 We declare that this manuscript is original, has not been published before and is not currently being considered for
publication elsewhere.

We confirm that the manuscript has been read and approved and that there are no other persons who satisfied the criteria for
authorship but are not listed.
 \section*{ Competing Interests:}

The authors declare that they have no known competing financial interests or personal relationships that could have appeared to influence the work reported in this paper.

\end{document}